# SPECTRAL CHARACTERIZATION OF AGING: THE REM-LIKE TRAP MODEL[1]

By Anton Bovier and Alessandra Faggionato

*Weierstrass Institut für Angewandte Analysis und Stochastik*

We review the aging phenomenon in the context of the simplest trap model, Bouchaud's REM-like trap model, from a spectral theoretic point of view. We show that the generator of the dynamics of this model can be diagonalized exactly. Using this result, we derive closed expressions for correlation functions in terms of complex contour integrals that permit an easy investigation into their large time asymptotics in the thermodynamic limit. We also give a "grand canonical" representation of the model in terms of the Markov process on a Poisson point process. In this context we analyze the dynamics on various time scales.

**1. Introduction.** The particular properties of the long term dynamics of many complex and/or disordered systems have been the subject of great interest in the physics and, increasingly, the mathematics community. The key paradigm here is the notion of *aging*, a notion that can be characterized in terms of scaling properties of suitable autocorrelation functions. Typically, aging can be associated to the existence of *infinitely many* time-scales that are inherently relevant to the system. In that respect, aging systems are distinct from *metastable* systems, which are characterized by a finite number of well separated time-scales, corresponding to the lifetimes of different metastable states.

Aging systems are rather difficult to analyze, both numerically and analytically. Most analytical results, even on the heuristic level, concern either the Langevin dynamics of spherical mean field spin glasses or *trap models*, a class of artificial Markov processes that in some way tries to mimic the long term dynamics of highly disordered systems (see, e.g., [8]).

Received May 2004; revised December 2004.
[1]Supported in part by the DFG in the Dutch-German Bilateral Research Group "Mathematics of Random Spatial Models from Physics and Biology."

One of the natural questions one is led to ask when being confronted with phenomena related to multiple time-scales is whether and how they can be related to *spectral properties*. This relationship has been widely investigated in the context of Markov processes with metastable behavior (see, e.g., [12, 13, 14, 20, 21, 10]), and it would be rather interesting to obtain a spectral characterization of aging systems as well, at least in the context of Markov processes. To our knowledge, this problem has not been widely studied so far. The only papers dealing with the problem are [24], by Butaud and Mélin, that have tackled one of the simplest trap models and on which we will comment below, and [17] and [23], that investigate convergence to equilibrium in the Random Energy Model (REM).

The present paper is intended to make a modest step in this direction by analyzing the relation between spectral properties and aging rigorously in the REM-like trap model. While this model may seem misleadingly simple, it has in the past provided valuable insights into the mechanisms of aging, and it is our hope that the analysis presented here will provide useful guidelines for further investigations of more complicated models.

The paper will be divided into two parts. In the first we analyze the REM-like trap model in the standard formulation of Bouchaud [9]. In the second part we go one step further and reformulate the model in a slightly different way as a Markov process on a Poisson point process. This formulation makes the relation to the real REM more suggestive (see [3, 4] for a full analysis), and allows, in a natural way, to study the dynamics of the model on different time scales.

**2. The REM-like trap model.** Let us recall the definition of *trap models* as introduced by Bouchaud and Dean [9]. Let $\mathcal{G} = (\mathcal{S}, \mathcal{E})$ be a finite graph with vertex set, $\mathcal{S}$, and edge set, $\mathcal{E}$. Let $\underline{E} := \{E_i, i \in \mathcal{S}\}$ be a random field, called *energy landscape*, and let $Y(t)$ be a continuous-time random walk on $\mathcal{G}$ with $\underline{E}$-dependent transition rates, $c_{i,j}$, such that $c_{i,j} > 0$ iff $\{i,j\} \in \mathcal{E}$, and

$$\mathbb{P}(Y(t+dt) = j | Y(t) = i) = c_{i,j} \, dt.$$

Setting $\tau_i^{-1} := \sum_{j \neq i} c_{i,j}$ and $p_{i,j} := c_{i,j} \tau_i$, the random walk, $Y(t)$, can be described as follows: after reaching the site $i$, the walk waits an exponential time of mean $\tau_i$ and then jumps to an adjacent site, $j$, with probability $p_{i,j}$. In the trap model, the transition rates are assumed to satisfy the following properties:

(2.1) $$e^{E_i} c_{i,j} = e^{E_j} c_{j,i} \qquad \forall \{i,j\} \in \mathcal{E},$$

(2.2) $$\mathbb{E}(\tau_i) = \infty,$$

where $\mathbb{E}$ denotes the expectation w.r.t. the random field $\underline{E}$. Since in several physical experiments (see [26]) the system is initially in equilibrium at a high temperature, $T \gg T_g$ and then is quickly cooled to a temperature smaller than $T_g$, and then its response to an external perturbation is measured, it is reasonable to consider $Y(t)$ with uniform initial distribution. A classical time–time correlation function is given by

$$\Pi(t, t_w) := \mathbb{P}(Y(s) = Y(t_w) \ \forall s \in [t_w, t_w + t]).$$

In order to observe aging, it is necessary to consider a thermodynamic limit, with the size of $\mathcal{G}$ going to infinity, and possibly a suitable time-rescaling. Rather recently, there have been a number of rigorous papers devoted to the analysis of trap models on the lattices $\mathbb{Z}$ [5, 18, 19] and $\mathbb{Z}^d$ [6, 11].

In this paper we consider the simplest trap model, called the *REM-like trap model* [9], that corresponds to choosing $\mathcal{G}$ to be the complete graph on $N$ vertices, that is,

$$\mathcal{G}_N = (\mathcal{S}_N, \mathcal{E}_N), \qquad \mathcal{S}_N := \{1, 2, \ldots, N\}, \qquad \mathcal{E}_N := \{\{i, j\} \ i \neq j \in \mathcal{S}_N\},$$

and to take as energy landscape a family, $\underline{E} = \{E_i : i \in \mathbb{N}\}$, of independent, exponentially distributed random variables, with parameter, $\alpha$, with $0 < \alpha < 1$. Given $N \in \mathbb{N}$, let $Y_N(t)$ be the continuous-time random walk on $\mathcal{G}_N$ with transition rates $c_{i,j} = e^{-E_i}/N$, for $i \neq j$. Setting $x_i = e^{-E_i}$, the infinitesimal generator of the random walk is given by

$$(2.3) \qquad \mathbb{L}_N := \begin{pmatrix} \frac{(N-1)x_1}{N} & -\frac{x_1}{N} & \cdots & -\frac{x_1}{N} \\ -\frac{x_2}{N} & \frac{(N-1)x_2}{N} & \cdots & -\frac{x_2}{N} \\ \vdots & \vdots & \ddots & \vdots \\ -\frac{x_N}{N} & -\frac{x_N}{N} & \cdots & \frac{(N-1)x_N}{N} \end{pmatrix}.$$

The dynamics can be described as follows: after reaching the state $i$, the walk waits an exponential time of mean $\frac{N}{N-1}e^{E_i}$ and then jumps with uniform probability to another state. Although, strictly speaking, the mean waiting time is given by $\frac{N}{N-1}e^{E_i}$, we call $\tau_i := x_i^{-1} = e^{E_i}$ waiting time (the discrepancy is negligible in the thermodynamic limit $N \uparrow \infty$).

Note that $\tau_i$ and $x_i$ have distributions given by

$$p(\tau) \, d\tau = \alpha \tau^{-1-\alpha} \, d\tau \qquad (\tau \geq 1); \qquad p(x) \, dx = \alpha x^{\alpha-1} \, dx \qquad (0 < x \leq 1),$$

respectively; in particular, $\mathbb{E}(\tau_i) = \infty$. Moreover, the equilibrium measure is given by $\mu_{eq}(i) = \tau_i / (\sum_{j=1}^N \tau_j)$. We are interested in the out-of-equilibrium dynamic with uniform initial distribution. $\mathbb{P}_N$ denotes the law of this random walk, given a realization of the random variables $E_i$.

Aging in the REM-like trap model is manifest from the asymptotic behavior of the time–time correlation function

$$(2.4) \qquad \Pi_N(t, t_w) := \mathbb{P}_N(Y_N(s) = Y_N(t_w) \ \forall s \in [t_w, t_w + t]).$$

Namely, as shown in [9], for almost all $\underline{E}$, and for all $\theta > 0$,

$$(2.5) \qquad \lim_{t_w \uparrow \infty} \lim_{N \uparrow \infty} \Pi_N(\theta t_w, t_w) = \frac{\sin(\pi\alpha)}{\pi} \int_{\theta/(1+\theta)}^1 u^{-\alpha}(1-u)^{\alpha-1}\,du.$$

Our main aim here is to show that the aging behavior of the system, derived in [9] using renewal arguments, can be obtained solely from *spectral information* about the generator $\mathbb{L}_N$. The method developed below will allow us to get further information on $Y_N(t)$ from the spectral properties of $\mathbb{L}_N$. In particular, given a function $h$ on $(0, \infty)$, it is possible to describe the asymptotic behavior of $\mathbb{E}_N(h(x_N(t)))$ and $\mathbb{E}_N(h(\tau_N(t)))$, where $\mathbb{E}_N$ denotes the expectation w.r.t. $\mathbb{P}_N$, and where $x_N(t), \tau_N(t)$ are defined as

$$x_N(t) = x_k, \qquad \tau_N(t) = \tau_k \qquad \text{if } Y_N(t) = k.$$

These results will allow us to investigate how the walk, as time goes on, visits deeper and deeper traps, that is, sites with larger and larger waiting time $\tau_i$ (see Section 2.2).

We start by giving a complete description of the eigenvalues and eigenvectors of $\mathbb{L}_N$. Let $\mu = \mu_N$ be the measure on $\mathcal{S}_N$ with $\mu(i) = x_i^{-1} = \tau_i$. Note that $\mathbb{L}_N$ is a symmetric operator on $L^2(\mu)$ and, trivially, $\mathbb{L}_N \mathbb{I} = 0$, where $\mathbb{I}$ is the vector with all entries equal to 1. The following proposition is based on elementary linear algebra:

PROPOSITION 2.1. *Let* $x_1, x_2, \ldots, x_N$ *be all distinct. Then,* $\mathbb{L}_N$ *has $N$ positive simple eigenvalues* $0 = \lambda_1 < \lambda_2 < \cdots < \lambda_N$ *such that*

$$\{\lambda_1, \lambda_2, \ldots, \lambda_N\} = \{\lambda \in \mathbb{C} : \phi(\lambda) = 0\},$$

*where $\phi(\lambda)$ is the meromorphic function*

$$(2.6) \qquad \phi(\lambda) := \sum_{j=1}^N \frac{\lambda}{x_j - \lambda}, \qquad (\lambda \in \mathbb{C}).$$

*If the $x_i$ are labelled such that $x_1 < x_2 < \cdots < x_N$, then $x_i < \lambda_{i+1} < x_{i+1}$, for $i = 1, \ldots, N-1$. Moreover, for any $i = 1, \ldots, N$, the vector $\psi^{(i)} \in \mathbb{R}^N$, defined as*

$$\psi_j^{(i)} := \frac{x_j}{x_j - \lambda_i} \qquad \text{for } j = 1, \ldots, N,$$

*is an eigenvector of $\mathbb{L}_N$ with eigenvalue $\lambda_i$. $\psi^{(1)}, \ldots, \psi^{(N)}$ form an orthogonal basis of $L^2(\mu)$.*

Since the $x_i$ have an absolutely continuous distribution, we trivially have the following:

COROLLARY 2.2. *The assertions of Proposition 2.1 hold with probability one for all $N$.*

PROOF. Let $\lambda$ be a generic eigenvalue and let us write the corresponding eigenvector, $\psi$, as $\psi = a(1,\ldots,1)^t + w$, where $\sum_{j=1}^{N} w_j = 0$. Since $(\mathbb{L}_N \psi)_j = x_j w_j$, we have to solve the system

$$x_j w_j = \lambda a + \lambda w_j \qquad \forall\, j = 1, \ldots, N. \tag{2.7}$$

Since $x_1, \ldots, x_N$ are distinct, it must be true that $a \neq 0$ (otherwise we get $\psi = 0$). Without loss of generality, we set $a = 1$. Note that $\lambda \neq x_j$, for $j = 1, \ldots, N$, since otherwise (2.7) would imply that $\lambda = 0 = x_j$. Therefore, we get $w_j = \frac{\lambda}{x_j - \lambda}$. Since it must be true that $\sum_{j=1}^{N} w_j = 0$, we get that $\lambda$ is an eigenvalue with $\psi$ s.t. $\psi_j = \frac{x_j}{x_j - \lambda}$ being the corresponding eigenvector, iff $\phi(\lambda) = 0$. This implies that $\phi$ has at most $N$ zeros. Since $\phi(0) = 0$, and, for real $\lambda$, $\lim_{\lambda \downarrow x_i} \phi(\lambda) = -\infty$, $\lim_{\lambda \uparrow x_i} \phi(\lambda) = \infty$, we get that $\phi$ has exactly $N$ zeros. From here the assertions of the theorem follow immediately. □

Proposition 2.1 has the following simple corollary:

COROLLARY 2.3. *With probability one, the spectral distribution $\sigma_N := \mathrm{Av}_{j=1}^{N} \delta_{\lambda_j}$ converges weakly to the measure $\alpha x^{\alpha - 1}\, dx$ on $[0,1]$.*

REMARK. The results of Proposition 2.1 are incompatible with the heuristic predictions in [24]. The discrepancy is particularly pronounced in the case of the eigenfunction. The reason for this is an inappropriate use of perturbation expansion in [24]. We will explain this in some detail in the Appendix.

We will now show that Proposition 2.1 allows to derive the asymptotics of the autocorrelation functions easily. In fact, it contains far more information on the long time behavior of the systems, some of which we will bring to light later.

Recall that $p_t(i,j)$, the probability to jump from $i$ to $j$ in an interval of time $t$, can be expressed as $p_t(i,j) = (e^{-t\mathbb{L}_N})_{i,j}$. In particular, by writing $\nu_t$ for the probability distribution of $Y_N(t)$ and thinking of the Radon derivative $\frac{d\nu_t}{d\mu}$ as column vector,

$$\frac{d\nu_t}{d\mu} = e^{-t\mathbb{L}_N} \frac{d\nu_0}{d\mu},$$

we see that

$$\frac{d\nu_t}{d\mu} = \sum_{k=1}^{N} \frac{\langle d\nu_0/d\mu, \psi^{(k)} \rangle}{\langle \psi^{(k)}, \psi^{(k)} \rangle} e^{-t\lambda_k} \psi^{(k)}. \tag{2.8}$$

The above formulas are true for an arbitrary initial distribution. Taking $\nu_0$ to be the uniform distribution, by Proposition 2.1, we get

$$\frac{d\nu_o}{d\mu} = \sum_{k=1}^N \gamma_k \psi^{(k)} \qquad \text{where } \gamma_k^{-1} := \langle \psi^{(k)}, \psi^{(k)} \rangle = \sum_{j=1}^N \frac{x_j}{(x_j - \lambda_k)^2}.$$

Then, by Proposition 2.1 and (2.8),

$$(2.9) \qquad \Pi_N(t, t_w) = \sum_{j=1}^N \sum_{k=1}^N \frac{\gamma_k e^{-\lambda_k t_w}}{x_j - \lambda_k} e^{-((N-1)/N)x_j t}$$

and

$$(2.10) \qquad \mathbb{E}_N(h(x_N(t))) = \sum_{j=1}^N \sum_{k=1}^N \frac{\gamma_k e^{-\lambda_k t}}{x_j - \lambda_k} h(x_j).$$

The above formulas (that may appear rather ugly at first sight) admit a nice complex integral representation through the following lemma:

LEMMA 2.4. *Let $\gamma$ be a positive oriented loop on $\mathbb{C}$ containing in its interior $\lambda_1, \ldots, \lambda_N$. Let $g$ be a holomorphic function on a domain $D \subset \mathbb{C}$ with $\gamma \subset D$. Then, for any $j = 1, \ldots, N$,*

$$(2.11) \qquad \sum_{k=1}^N \frac{\gamma_k g(\lambda_k)}{x_j - \lambda_k} = \frac{1}{2\pi i} \int_\gamma \frac{g(\lambda)}{\phi(\lambda)(x_j - \lambda)} d\lambda.$$

PROOF. Let us set $X := \{x_1, \ldots, x_N\}$ and $\Lambda := \{\lambda_1, \lambda_2, \ldots, \lambda_N\}$. Then, $\phi(\lambda)$ is a holomorphic function on $\mathbb{C} \setminus X$, where $\phi'(\lambda) = \sum_{j=1}^N \frac{x_j}{(x_j - \lambda)^2}$, and, in particular, $\phi'(\lambda_j) = \gamma_j^{-1}$. Moreover, the function $[\phi(\lambda)(x_j - \lambda)]^{-1}$, a priori defined on $\mathbb{C} \setminus (X \cup \Lambda)$, can be analytically continued to $X$ as a meromorphic function with simple poles only at the points of $\Lambda$. Now the conclusion follows from a trivial application of the residue theorem. □

We can obviously use Lemma 2.4 to rewrite (2.9) and (2.10) in the form

$$(2.12) \quad \Pi_N(t, t_w) = \frac{1}{2\pi i} \int_\gamma \frac{e^{-t_w \lambda}}{\lambda} \left( \text{Av}_j \frac{e^{-(N-1)/N x_j t}}{x_j - \lambda} \bigg/ \text{Av}_j \frac{1}{x_j - \lambda} \right) d\lambda,$$

$$(2.13) \quad \mathbb{E}_N(h(x_t)) = \frac{1}{2\pi i} \int_\gamma \frac{e^{-t\lambda}}{\lambda} \left( \text{Av}_j \frac{h(x_j)}{x_j - \lambda} \bigg/ \text{Av}_j \frac{1}{x_j - \lambda} \right) d\lambda,$$

where $\text{Av}_j$ denotes the average over $j = 1, 2, \ldots, N$.

The above integral representations of $\Pi_N(t, t_w)$ and $\mathbb{E}_N(h(x_t))$ have two advantages. First, the appearance of averages allows to compute their limiting behavior as $N \uparrow \infty$ easily by using the ergodicity of the random field $\underline{E}$. Second, by means of the residue theorem, their Laplace transform can be easily computed in order to derive the asymptotic behavior of $\Pi_N(t, t_w)$ and $\mathbb{E}_N(h(x_t))$, for $N, t_w, t \gg 1$ (see Sections 2.1 and 2.2).

A much more general derivation of the above integral representations is discussed in Appendix A.3.

### 2.1. Aging behavior of $\Pi_N(t, t_w)$.

PROPOSITION 2.5. *Let us define*

$$(2.14) \qquad \Pi(t, t_w) := \frac{1}{2\pi i} \int_\gamma \frac{e^{-t_w \lambda}}{\lambda} \frac{\mathbb{E}_x(e^{-xt}/(\lambda - x))}{\mathbb{E}_x(1/(\lambda - x))} \, d\lambda,$$

*where $\mathbb{E}_x$ is the expectation w.r.t. the measure $\alpha x^{\alpha-1} \, dx$ on $[0, 1]$ and $\gamma$ is any positive oriented complex loop around the interval $[0, 1]$. Then,*

$$(2.15) \qquad \lim_{N \uparrow \infty} \Pi_N(t, t_w) = \Pi(t, t_w) \qquad \forall t, t_w, \ a.s.$$

PROOF. Recall (2.12) and fix $0 < \delta < 1/2$. Due to analyticity, we can choose the integration contour, $\gamma$, to have distance 1 from the segment $[0, 1]$. For each $\lambda \in \gamma$, the random variables $(x_j - \lambda)^{-1}$, $j \in \mathbb{N}$, are i.i.d. and bounded. Therefore, for a suitable positive constant $c > 0$,

$$(2.16) \quad \mathbb{P}\left(\left|\mathrm{Av}_{j=1}^N \frac{1}{x_j - \lambda} - \mathbb{E}_x\left(\frac{1}{\lambda - x}\right)\right| \geq N^{-1/2+\delta}\right) \leq e^{-cN^{2\delta}} \qquad \forall \lambda \in \gamma.$$

Since for each $x \in [0, 1]$ and $\lambda \in \gamma$, $|\frac{\partial}{\partial \lambda}(x - \lambda)^{-1}| \leq 1$, a simple chaining argument allows to deduce from the pointwise estimate (2.16) uniform control in $\lambda$. Using the Borel–Cantelli lemma, one can then infer that, a.s.,

$$(2.17) \quad \sup_{\lambda \in \gamma} \left|\mathrm{Av}_{j=1}^N \frac{1}{x_j - \lambda} - \mathbb{E}_x\left(\frac{1}{\lambda - x}\right)\right| \leq cN^{-1/2+\delta} \qquad \forall N \in \mathbb{N}.$$

Similar arguments show that, a.s., given $M \in \mathbb{N}$, there exists a constant, $c_M$, such that

$$(2.18) \quad \sup_{M^{-1} \leq t \leq M} \sup_{\lambda \in \gamma} \left|\mathrm{Av}_{j=1}^N \frac{e^{-((N-1)/N)x_j t}}{x_j - \lambda} - \mathbb{E}_x\left(\frac{e^{-x_j t}}{\lambda - x_j}\right)\right| \leq c_M N^{-1/2+\delta}$$

$$\forall N \in \mathbb{N}.$$

Note that, for each $\lambda \in \gamma$, $\mathrm{Av}_{j=1}^N (x_j - \lambda)^{-1}$ is a convex combination of points of modulus larger or equal than $1/2$, contained in an angular sector with angle non-larger than a suitable constant, $c < \pi$. In particular, $|\mathrm{Av}_{j=1}^N (x_j - \lambda)^{-1}| \geq c' > 0$, for all $N$. From here the assertion of the proposition follows from Lebesgue's dominated convergence theorem. $\square$

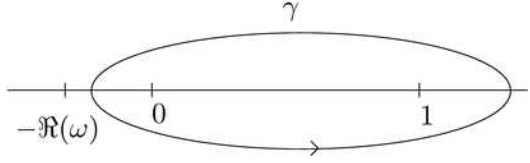

Fig. 1. *Possible integration contour.*

Given $\theta > 0$, we are interested in the limit of $\Pi(\theta t_w, t_w)$, as $t_w \uparrow \infty$. This will be done using the Laplace transform of $\Pi(\theta t_w, t_w)$,

$$\hat{\Pi}(\theta, \omega) := \int_0^\infty e^{-\omega t_w} \Pi(\theta t_w, t_w) \, dt_w, \qquad [\Re(\omega) > 0].$$

The computation of this Laplace transform is trivial, if we use the integral expression (2.14).

Let $\omega \in \mathbb{C}$, with $\Re(\omega) > 0$, and fix a positive oriented loop, $\gamma$, around the segment $[0, 1]$, such that $\gamma \subset \{z \in \mathbb{C} : \Re(z) > -\Re(\omega)\}$; see Figure 1. Then, $\Re(\omega + \lambda + x\theta) > 0$, for $x \in [0, 1]$ and $\lambda \in \gamma$, so that (2.14) implies

$$\hat{\Pi}(\theta, \omega) = \mathbb{E}_x \bigg( \frac{1}{2\pi i} \int_\gamma \bigg[ \lambda(\lambda - x)(\lambda + \omega + \theta x) \mathbb{E}_{\bar{x}} \bigg( \frac{1}{\lambda - \bar{x}} \bigg) \bigg]^{-1} d\lambda \bigg).$$

Here $\mathbb{E}_x$ and $\mathbb{E}_{\bar{x}}$ denote the expectation w.r.t. the measure $\alpha x^{\alpha-1} \, dx$ on $[0, 1]$.

Let us consider the change of variables $z = \frac{1}{\lambda}$ and write $\hat{\gamma}$ for the path $\gamma$ with inverted orientation (i.e., positive oriented w.r.t. $\lambda = \infty$). Then we get

$$\hat{\Pi}(\theta, \omega) = \mathbb{E}_x \bigg( \frac{1}{2\pi i} \int_{\hat{\gamma}} \bigg[ (1 - zx)(1 + z\omega + z\theta x) \mathbb{E}_{\bar{x}} \bigg( \frac{1}{1 - z\bar{x}} \bigg) \bigg]^{-1} dz \bigg).$$

Given $x \in [0, 1]$, the integrand is a meromorphic function in $\mathbb{C} \setminus [1, \infty)$ that has only a single pole of order 1 inside $\hat{\gamma}$, namely, at $z = -(\omega + x\theta)^{-1}$. By the residue theorem, we get

$$(2.19) \qquad \hat{\Pi}(\theta, \omega) = \mathbb{E}_x \bigg( \frac{1}{\omega + x\theta + x} \bigg/ \mathbb{E}_{\bar{x}} \bigg( \frac{\omega + x\theta}{\omega + x\theta + \bar{x}} \bigg) \bigg).$$

LEMMA 2.6. *The r.h.s. of* (2.19) *is well defined and holomorphic for any* $\omega \in \mathbb{C} \setminus (-\infty, 0]$. *In particular, the function* $\hat{\Pi}(\theta, \omega)$, *defined for* $\Re(\omega) > 0$, *can be analytically continued to the set* $\mathbb{C} \setminus (-\infty, 0]$.

PROOF. As proved in [15], Chapter 3, the Laplace transform, $\hat{\Pi}(\theta, \omega)$, is holomorphic on the set of convergence points. Therefore, we only need to show that the r.h.s. of (2.19) is well defined and holomorphic on $\Im(\omega) \neq 0$. Let us assume that $\Im(\omega) > a > 0$. Then, trivially, $\forall \, x, \bar{x} \in [0, 1]$,

$$\frac{\omega + x\theta}{\omega + x\theta + \bar{x}} \in \mathcal{B} := \{z \in \mathbb{C} : z = |z|e^{i\theta} \text{ with } 0 \leq \theta \leq \theta_0, |z| \geq c\},$$

SPECTRAL CHARACTERIZATION OF AGING 9for suitable constants, $c, \theta_0$, depending on $a$ and satisfying $\theta_0 < \pi$. Moreover, since $\lim_{|\omega|\uparrow\infty} \frac{\omega + x\theta}{\omega + x\theta + \bar{x}} = 1$,

$$(2.20) \quad 0 < c_1(a) \leq \left|\frac{\omega + x\theta}{\omega + x\theta + \bar{x}}\right| \leq c_2(a) \qquad \forall a > 0, \forall \omega : \Im(\omega) \geq a.$$

By (2.20) and the geometry of $\mathcal{B}$, we have that $\mathbb{E}_{\bar{x}}(\frac{\omega + x\theta}{\omega + x\theta + \bar{x}})$ is well defined and has distance $c_3(a) > 0$ from the origin. Moreover, $|\omega + x\theta + x| \geq \Im(\omega)$. Therefore, the r.h.s. of (2.19) is well defined and, due to the previous estimates and Lebesgue's dominated convergence theorem, it is continuous on $\{\Im(\omega) \neq 0\}$, thus implying continuity on $\mathbb{C} \setminus (-\infty, 0]$.

We recall Morera's theorem: if $f(\omega)$ is defined and continuous in a open set $\Omega \subset \mathbb{C}$ and if $\int_\gamma f\, d\omega = 0$, for all closed curves, $\gamma$, in $\Omega$, then $f(\omega)$ is holomorphic in $\Omega$. Therefore, using Fubini's and Morera's theorems, one can prove that the function $\mathbb{E}_{\bar{x}}(\frac{\omega + x\theta}{\omega + x\theta + \bar{x}})$ is holomorphic on $\mathbb{C} \setminus (-\infty, 0]$. The proof can be concluded by a second application of the same theorems. □

In what follows, we keep the notation, $\hat{\Pi}(\theta, \omega)$, for the analytic continuation of the Laplace transform. The next lemma describes the behavior of $\hat{\Pi}(\theta, \omega)$ near the origin. Using the Laplace inversion formula, we then derive from this result the asymptotic behavior of $\Pi(\theta, t_\omega)$, as $t_\omega \uparrow \infty$.

LEMMA 2.7. *For any $\theta > 0$, set*

$$A(\theta) =: \frac{\sin(\pi\alpha)}{\pi} \int_{\theta/(\theta+1)}^1 u^{-\alpha}(1-u)^{\alpha-1}\, du.$$

*Moreover, define*

$$(2.21) \qquad \mathcal{A} := \{re^{i\phi} : r \geq 0, |\phi| \leq \tfrac{3}{4}\pi\}.$$

*Then, for a suitable positive constant, $c > 0$,*

$$(2.22) \qquad |\hat{\Pi}(\theta, \omega)| \leq c|\omega|^{-1} \qquad \forall \omega \in \mathcal{A} : |\omega| \geq 1,$$

$$(2.23) \qquad |\hat{\Pi}(\theta, \omega) - A(\theta)/\omega| \leq c|\omega|^{-\alpha} \qquad \forall \omega \in \mathcal{A} : |\omega| \leq 1.$$

PROOF. The first estimate (2.22) follows trivially from (2.19) and (2.20). Let us prove (2.23) for $\omega \in \mathcal{A}$ and $|\omega| \leq 1$.

In what follows, $c_0, c_1, \ldots$ denote positive constants depending only on $\theta$. Moreover, given $z \in \mathbb{C}$, we denote by $\int_0^z$ and $\int_z^\infty$ the integrals over the paths $\{sz : 0 \leq s \leq 1\}$ and $\{sz : s \geq 1\}$, respectively. We extend the functions $z^{-\alpha}$ and $z^{\alpha-1}$, defined on $(0, \infty)$, to $\mathbb{C} \setminus (-\infty, 0]$ by analytic continuation. Then, (2.19) implies

$$\omega\hat{\Pi}(\theta, \omega)$$

$$(2.24) \quad = \int_0^{1/\omega} x^{\alpha-1} \left( [1+x(1+\theta)][1+x\theta] \int_0^{1/\omega} \frac{y^{\alpha-1}}{1+x\theta+y} dy \right)^{-1} dx$$

$$= \int_0^{1/\omega} x^{\alpha-1} \left( [1+x(1+\theta)][1+x\theta]^\alpha \int_0^{1/\omega(1+x\theta)} \frac{y^{\alpha-1}}{1+y} dy \right)^{-1} dx.$$

Define

$$\mathcal{B} := \left\{ (\omega, x) \in \mathbb{C}^2 \text{ s.t. } \omega \in \mathcal{A}, x = \frac{s}{\omega} \text{ for some } s : 0 \le s \le 1 \right\}.$$

Since $(\omega(1+x\theta))^{-1} \in \mathcal{A} \cap \{z : |z| \ge c_0\}$, we obtain

$$(2.25) \quad \left| \int_0^{1/(\omega(1+x\theta))} \frac{y^{\alpha-1}}{1+y} dy \right| \ge c_1,$$

$$(2.26) \quad \left| \int_{1/\omega(1+x\theta)}^\infty \frac{y^{\alpha-1}}{1+y} dy \right| \le c_2 |\omega(1+x\theta)|^{1-\alpha}.$$

Let $B(\alpha)$ be defined as

$$(2.27) \quad B(\alpha) := \int_0^\infty \frac{y^{\alpha-1}}{1+y} dy = \int_0^1 u^{-\alpha}(1-u)^{\alpha-1} du = \frac{\pi}{\sin(\pi\alpha)}$$

[note that the above second identity follows from the change of variables $u = y(1+y)^{-1}$, while the last one is well known in the theory of the Gamma function]. Using (2.25) and (2.26), we obtain

$$(2.28) \quad \begin{aligned} & \left| \omega \hat{\Pi}(\theta, \omega) - \frac{1}{B(\alpha)} \int_0^{1/\omega} \frac{x^{\alpha-1} dx}{(1+x(1+\theta))(1+x\theta)^\alpha} \right| \\ & \le |\omega|^{1-\alpha} \int_0^{1/\omega} \frac{|x|^{\alpha-1} d|x|}{|1+x(1+\theta)||1+x\theta|^{2\alpha-1}} \le c_3 |\omega|^{1-\alpha}. \end{aligned}$$

Since

$$\left| \int_{1/\omega}^\infty \frac{x^{\alpha-1} dx}{(1+x(1+\theta))(1+x\theta)^\alpha} \right| \le c_4 |\omega|$$

and, using analyticity and integrability of the singularities around $z=0$ and $z=\infty$,

$$\int_{s\omega : s \ge 0} \frac{x^{\alpha-1} dx}{(1+x(1+\theta))(1+x\theta)^\alpha} = \int_0^\infty \frac{x^{\alpha-1} dx}{(1+x(1+\theta))(1+x\theta)^\alpha},$$

we get

$$\left| \omega \hat{\Pi}(\theta, \omega) - \frac{1}{B(\alpha)} \int_0^\infty \frac{x^{\alpha-1}}{(1+x(1+\theta))(1+x\theta)^\alpha} dx \right| \le c_5 |\omega|^{1-\alpha}.$$

Using the change of variables $v = x^{-1} + \theta$ and $u = v(1+v)^{-1}$, we obtain

$$\int_0^\infty \frac{x^{\alpha-1}}{(1+x(1+\theta))(1+x\theta)^\alpha}\,dx = \int_{\theta/(\theta+1)}^1 u^{-\alpha}(1-u)^{\alpha-1}\,du,$$

which implies the assertion of the lemma. □

Lemma 2.7 and Proposition A.1 allow us to conclude the proof of the aging behavior of $\Pi_N(t, t_w)$, that is, to recover the result of Bouchaud and Dean [9]:

PROPOSITION 2.8. *For almost all energy landscapes $\underline{E}$, given $\theta > 0$,*

$$(2.29) \quad \lim_{t_w \uparrow \infty} \lim_{N \uparrow \infty} \Pi_N(\theta t_w, t_w) = \frac{\sin(\pi\alpha)}{\pi} \int_{\theta/(1+\theta)}^1 u^{-\alpha}(1-u)^{\alpha-1}\,du.$$

2.2. *Visiting deeper and deeper traps.* In this section we use the integral representation (2.13) in order to study the probability that the walk at time $t$ is in a deep trap, that is, in a state with large waiting time. In Proposition 2.9 we first prove that the probability to be in a site with waiting time smaller than $O(1)$ decays as $t^{\alpha-1}$, thus implying the aging behavior of correlation functions described in Section 3. In the second part, we will investigate the random variable $tx_N(t)$ and show that, for almost all $\underline{E}$, it has a weak limit, as first $N \uparrow \infty$ and then $t \uparrow \infty$. As consequence, with high probability, at time $t$, the system is in a state of waiting time $O(t)$, as stated in Proposition 2.10.

Reasoning as in the proof of Proposition 2.5, we can prove, for almost all energy landscapes, $\underline{E}$, that, given a function $h$ on $[0, 1]$, that can be uniformly approximated by piecewise $C^1$ functions,

$$(2.30) \quad \begin{aligned} H(t) &:= \lim_{N \uparrow \infty} \mathbb{E}_N(h(x_N(t))) \\ &= \frac{1}{2\pi i} \int_\gamma \frac{e^{-t\lambda}}{\lambda} \frac{\int_0^1 (h(x)/(\lambda-x))x^{\alpha-1}\,dx}{\int_0^1 (1/(\lambda-x))x^{\alpha-1}\,dx}\,d\lambda \qquad \forall t > 0, \end{aligned}$$

where $\gamma$ is a positive oriented loop around $[0, 1]$.

Since $H(t)$ is a bounded function, the Laplace integral $\hat{H}(\omega) := \int_0^\infty H(t)e^{-\omega t}\,dt$ is absolutely convergent when $\Re(\omega) > 0$. By the same arguments we used to derive (2.19), it is easy to deduce from the integral representation (2.30) that

$$(2.31) \quad \hat{H}(\omega) = \frac{1}{\omega} \frac{\int_0^1 (h(x)/(\omega+x))x^{\alpha-1}\,dx}{\int_0^1 (1/(\omega+x))x^{\alpha-1}\,dx}.$$

In the following proposition we concentrate on the case $h(x) := \mathbb{I}_{x \geq \delta}$. By (2.31), we can give precise information on the asymptotic behavior of the probability to be at time $t$ in a site with waiting time smaller than $1/\delta$:

PROPOSITION 2.9. *Let*

$$B(\delta) := \frac{\int_\delta^1 x^{\alpha-2}\,dx}{\int_0^\infty x^{\alpha-1}/(1+x)\,dx}, \qquad c(\alpha) := \int_0^\infty y^{\alpha-1} e^{-y}\,dy. \tag{2.32}$$

*Then, for almost all energy landscapes,* $\underline{E}$,

$$\lim_{s\uparrow\infty} s^{1-\alpha} \lim_{N\uparrow\infty} \mathbb{P}_N(x_N(s) > \delta) = B(\delta)/c(\alpha). \tag{2.33}$$

Finally, we show that, with high probability, at time $t$, the walk is in a trap of depth of order $O(t)$. In particular, the random variables $tx_N(t)$ converge weakly to a nonnegative random variable, as first $N \uparrow \infty$ and then $t \uparrow \infty$, a.s. This result corresponds to the convergence of the expectation of bounded, continuous functions and, due to Lemma 2.11, such a convergence can be extended to the larger class of bounded, piecewise continuous functions, which is more suitable for the investigation of the phenomenon of visiting deeper and deeper traps:

PROPOSITION 2.10. *Let $Z$ be the unique random variable with range in $(0,\infty)$ having Laplace transform*

$$\mathbb{E}(e^{-\theta Z}) = \frac{\sin(\pi\alpha)}{\pi} \int_{\theta/(\theta+1)}^1 u^{-\alpha}(1-u)^{\alpha-1}\,du.$$

*Then, for almost all energy landscape $\underline{E}$, given a bounded piecewise continuous function, $h$, on $(0,\infty)$,*

$$\lim_{t\uparrow\infty} \lim_{N\uparrow\infty} \mathbb{E}_N(h(tx_N(t)))$$

$$= \lim_{t\uparrow\infty} \frac{1}{2\pi i} \int_\gamma \frac{e^{-t\lambda}}{\lambda} \frac{\int_0^1 (h(xt)/(\lambda-x)) x^{\alpha-1}\,dx}{\int_0^1 (1/(\lambda-x)) x^{\alpha-1}\,dx}\,d\lambda \tag{2.34}$$

$$= \mathbb{E}(h(Z)).$$

*In particular, for almost all energy landscapes,* $\underline{E}$,

$$\lim_{t\uparrow\infty} \lim_{N\uparrow\infty} \mathbb{P}\left(\frac{\tau_N(t)}{t} \geq u\right) = \mathbb{P}(Z \leq u^{-1}) \qquad \forall u > 0.$$

PROOF OF PROPOSITION 2.9. We have to prove that $\lim_{s\uparrow\infty} s^{1-\alpha} H(s) = B(\delta)/c(\alpha)$, where $H$ is given by (2.30) with $h(x) := \mathbb{I}_{x \geq \delta}$. As in the proof of Lemma 2.6, we can show that the r.h.s. of (2.31) is well defined and holomorphic on $\mathbb{C} \setminus (-\infty, 0]$. We keep the notation $\hat{H}$ for this extended function. Changing variables $x = \omega y$, we get

$$\int_0^1 \frac{x^{\alpha-1}}{\omega+x}\,dx = \omega^{\alpha-1} \int_{\gamma_w} \frac{y^{\alpha-1}}{1+y}\,dy, \tag{2.35}$$

where $\gamma_\omega$ is the oriented path $\{s/\omega\}_{0 \le s \le 1}$. Let $\hat\gamma_\omega$ be the path $\{s/\omega\}_{s \ge 0}$. By analyticity and integrability of the singularities at $z = 0, z = \infty$, we have

$$\int_{\hat\gamma_\omega} \frac{y^{\alpha-1}}{1+y}\,dy = \int_0^\infty \frac{y^{\alpha-1}}{1+y}\,dy.$$

Let us define $\mathcal{A} := \{re^{i\theta} : 0 < r < \infty, |\theta| \le \tfrac{3}{4}\pi\}$. Then, for a suitable constant $c_1$,

$$\left|\int_{\hat\gamma_\omega \setminus \gamma_\omega} \frac{y^{\alpha-1}}{1+y}\,dy\right| \le c_1 |\omega|^{1-\alpha} \qquad \forall\, \omega \in \mathcal{A} : |\omega| \le 1,$$

implying

(2.36) $$\int_0^1 \frac{x^{\alpha-1}}{\omega + x}\,dx = \omega^{\alpha-1}\left(\int_0^\infty \frac{y^{\alpha-1}}{1+y}\,dy + O(|\omega|^{1-\alpha})\right),$$

where $A = B + O(1/N)$ is understood to mean that there exists $C < \infty$ such that $|A - B| \le C/N$. Trivially,

(2.37) $$\int_\delta^1 \frac{x^{\alpha-1}}{\omega + x}\,dx = (1 + O(|\omega|))\int_\delta^1 x^{\alpha-2}\,dx.$$

Note that the estimate of the error terms in (2.36) and (2.37) is uniform in $\omega \in \mathcal{A}$, $|\omega| \le 1$. Then, from (2.31), (2.36) and (2.37), we get

(2.38) $$|\omega^\alpha \hat{H}(\omega) - B(\delta)| \le c_2 |\omega|^{1-\alpha} \qquad \forall\, \omega \in \mathcal{A} : |\omega| \le 1.$$

Since, trivially, $|\hat{H}(\omega)| \le c_3 |\omega|^{-1}$, for $\omega \in \mathcal{A}$ with $|\omega| > 1$, the assertion of the proposition follows from Proposition A.1. $\square$

PROOF OF PROPOSITION 2.10. As discussed before (2.30), one can show that, for almost all energy landscapes, $\underline{E}$, given a piecewise continuous function $h$ on $(0, \infty)$,

$$\Phi_t(h) := \lim_{N\uparrow\infty} \mathbb{E}_N(h(tx_N(t)))$$

$$= \frac{1}{2\pi i}\int_\gamma \frac{e^{-t\lambda}}{\lambda} \frac{\int_0^1 (h(xt)/(\lambda - x))x^{\alpha-1}\,dx}{\int_0^1 (1/(\lambda - x))x^{\alpha-1}\,dx}\,d\lambda \qquad \forall\, t \ge 0,$$

where $\gamma$ is a positive oriented loop around $[0, 1]$. Note that $\Phi_t$ defines a positive linear functional on the space of continuous functions on $(0, \infty)$ that decay at $\infty$ and satisfy $\Phi_t(1) = 1$. Therefore, the Riesz–Markov representation theorem (see Theorem IV.18 in [25]) implies that $\Phi_t(h) = \mu_t(h)$, for a unique Borel probability measure $\mu_t$ on $[0, \infty)$. In particular, there exists a random variable, $Z_t$ on $(0, \infty)$, such that

$$\lim_{N\uparrow\infty} tx_N(t) \to Z_t \qquad \text{weakly, } \forall\, t > 0 \text{ a.s.}$$

If we take $h(t) = e^{-t\theta}$, then $\Phi_t(h) = \mu_t(h) = \Pi(\theta t, t)$, with $\Pi$ defined as in (2.14). That means that $\Pi(\theta t, t)$ is the Laplace transform of $Z_t$. As proved in Section 2.1,

$$\lim_{t\uparrow\infty} \Pi(\theta t, t) = \frac{\sin(\pi\alpha)}{\pi} \int_{\theta/(\theta+1)}^{1} u^{-\alpha}(1-u)^{\alpha-1}\,du := f(\theta).$$

We claim that $f(\theta)$ is the Laplace transform of a random variable $Z$ with range in $[0,\infty)$. To show this, we apply the criterion given by Theorem 1, Section XIII.4 in [16]. By (2.27), $f(0) = 1$. Moreover, $f^{(1)}(\theta) = -\frac{\sin(\pi\alpha)}{\pi}\theta^{-\alpha}(1+\theta)^{-1}$, thus implying (by trivial computations) that $(-1)^k f^{(k)}(\theta) \geq 0$. This completes the proof of our statement.

Since the Laplace transform of $Z_t$ converges to the Laplace transform of $Z$, as $t \uparrow \infty$, it follows that $Z_t$ converges weakly to $Z$, implying (2.34), whenever $h$ is a bounded continuous functions on $(0,\infty)$. Finally, due to Theorem 5.2 in [7], convergence still holds if $h$ is a bounded measurable function, whose set of discontinuity points has zero measure w.r.t. the distribution of $Z$. Therefore, Lemma 2.11 allows to prove (2.34) for $h$ bounded and piecewise continuous. $\square$

LEMMA 2.11. *The distribution function, $F(z) := \mathbb{P}(Z \leq z)$, of the positive random variable $Z$ is continuous.*

PROOF. Trivially, $F$ is increasing and right continuous. Therefore, it has a countable set of points of discontinuity. Moreover, by the Laplace inversion formula (see [16], XIII.4), if $x$ is a point of continuity, then

$$F(x) = \lim_{a\to\infty} \sum_{n\leq ax} \frac{(-a)^n}{n!} f^{(n)}(a).$$

Given $s = 0, 1, 2, \ldots$ and $\gamma > 0$, let $c_s(\gamma) > 0$ be such that $D_a^s a^{-\gamma} = (-1)^s \times c_s(\gamma) a^{-\gamma-s}$. Then the Leibniz formula implies

$$(-1)^n D_a^n (a^{-\alpha}(1+a)^{-1}) = \sum_{s=0}^{n} c_s(\alpha) c_{n-s}(1) a^{-\alpha-s}(1+a)^{-1-n+s}$$

$$\leq (-1)^n D_a^n a^{-\alpha-1}.$$

Since

$$f^{(1)}(a) = -\frac{\sin(\pi\alpha)}{\pi} a^{-\alpha}(1+a)^{-1},$$

the above estimate implies

$$(-1)^n f^{(n)}(a) = |f^{(n)}(a)| \leq \frac{\sin(\pi\alpha)}{\pi} \prod_{k=1}^{n-1}(k+\alpha) a^{-n-\alpha} \qquad \forall\, n \geq 2.$$

In particular, given two points of continuity, $0 < x < z$, we have

$$(2.39) \quad F(z) - F(x) \leq \frac{\sin(\pi\alpha)}{\pi} \limsup_{a \to \infty} a^{-\alpha} \sum_{ax < n \leq az} \frac{1}{n} \prod_{k=1}^{n-1} \left(1 + \frac{\alpha}{k}\right).$$

One can prove that the sequence $\prod_{k=1}^{n-1} e^{-\alpha/n}(1 + \frac{\alpha}{k})$ is convergent (see [1], Chapter 5, Section 2.4). Denote its limit by $c_\alpha$ and let $\gamma$ be Euler's constant

$$\gamma := \lim_{n \uparrow \infty} \left(1 + \frac{1}{2} + \frac{1}{3} + \cdots + \frac{1}{n} - \log n\right).$$

Then we can write

$$(2.40) \quad \frac{1}{n} \prod_{k=1}^{n-1}\left(1 + \frac{\alpha}{k}\right) = e^{\alpha(1+1/2+\cdots+1/(n-1)-\log(n-1))} \frac{(n-1)^\alpha}{n} \prod_{k=1}^{n-1} e^{-\alpha/k}\left(1 + \frac{\alpha}{k}\right).$$

In particular, we can substitute in (2.40) $\prod_{k=1}^{n-1} e^{-\alpha/k}(1 + \frac{\alpha}{k})$ with $c_\alpha$ with an error term in (2.39) bounded by

$$\text{const.} a^{-\alpha}(az - ax)\frac{(az)^\alpha}{ax}\left|\prod_{k=1}^{n-1} e^{-\alpha/k}\left(1 + \frac{\alpha}{k}\right) - c_\alpha\right|$$

$$\leq c(x,z) \left|\prod_{k=1}^{n-1} e^{-\alpha/k}\left(1 + \frac{\alpha}{k}\right) - c_\alpha\right|,$$

which is negligible, as $a \uparrow \infty$. Therefore,

$$(2.41) \quad F(z) - F(x) \leq \text{const.} \limsup_{a \to \infty} c_\alpha a^{-\alpha} e^{\gamma\alpha} \sum_{ax < n \leq az} (n-1)^{-1+\alpha}$$
$$\leq c'(z^\alpha - x^\alpha),$$

for some positive constant $c'$. Since (2.41) is valid almost everywhere and $F$ is monotone, it follows that $F$ is continuous. $\square$

**3. Other correlation functions.** In this section we study the asymptotic behavior of different time–time correlation functions, $\Pi_N^{(1)}(t, t_w)$, $\Pi_N^{(2)}(t, t_w)$, for which *deep* traps play a special role. This section is mainly a preparation of what is to follow in the second part of the paper.

Given $\delta > 0$, we define the set of sites with small waiting time as $D_N := \{i : x_i \geq \delta, i = 1, \ldots, N\}$. Moreover, we set

$$(3.1) \quad \Pi_N^{(1)}(t, t_w) := \mathbb{P}_N(Y_N(u) \in D_N \ \forall u \in (t_w, t_w + t] \text{ s.t. } Y_N(u) \neq Y_N(u^-))$$

$$(3.2) \quad \Pi_N^{(2)}(t, t_w) := \mathbb{P}_N(Y_N(u) \in D_N \cup \{Y_N(t_w)\}$$
$$\forall u \in (t_w, t_w + t] \text{ s.t. } x_N(u) \neq x_N(u^-)).$$

Given a subset $A \subset \mathcal{S}_N$, $i \in A$ and $s > 0$, let $\varphi_{N,A}(i,s)$ be defined as

$$\varphi_{N,A}(i,s) := \mathbb{P}_N(Y_N(u) \in A \ \forall u \in [0,s] | Y_N(0) = i).$$

Then

$$\Pi_N^{(1)}(t,t_w) = \Pi_N(t,t_w)$$

(3.3)
$$+ \sum_{j=1}^N \mathbb{P}_N(Y_N(t_w) = j) \int_0^t ds \frac{x_j e^{-sx_j}}{N} \sum_{i \in D_N} \varphi_{N,D_N}(i, t-s),$$

(3.4) $\Pi_N^{(2)}(t,t_w) = \sum_{j=1}^N \mathbb{P}_N(x_N(t_w) = x_j) \varphi_{N, D_N \cup \{j\}}(j,t),$

where the first identity can be derived by conditioning on the first jump performed after the waiting time $t_w$ and by recalling the following realization of the dynamics: after arriving at the state $i$, the system waits an exponential time with parameter $x_i$ and after that it jumps to a site in $\mathcal{S}_N$ with uniform probability.

The following proposition is mainly a consequence of the phenomenon of visiting deep traps with higher and higher probability. Recall that Proposition 2.10 implies

(3.5)  $\qquad \lim_{t \uparrow \infty} \lim_{N \uparrow \infty} \mathbb{P}_N(x_N(t) > \varepsilon) = 0 \qquad \forall \varepsilon > 0.$

PROPOSITION 3.1.  *For almost all $\underline{x}$,*

(3.6) $\qquad \lim_{t_w \uparrow \infty} \sup_{t \geq 0} |\Pi_N^{(i)}(t,t_w) - \Pi_N(t,t_w)| = 0 \qquad \textit{for } i = 1,2.$

PROOF. We consider first the case $i = 1$.
We claim that, for any $u > 0$ and $i \in D_N$,

(3.7) $\qquad \varphi_{N,D_N}(i,u) \leq \exp\left(-\delta u \left(1 - \frac{|D_N|}{N}\right)\right).$

In order to prove such a bound, we introduce a new random walk, $Y_N^*(t)$, whose generator, $\mathbb{L}^*$, defined as the r.h.s. of (2.3) with $x_i$ replaced by $\delta$ if $i \in D_N$. By a simple coupling argument, one gets

$$\varphi_{N,D_N}(i,u) \leq \varphi_{N,D_N}^*(i,u),$$

where the function $\varphi_{N,D_N}^*$ is the analogue of $\varphi_{N,D_N}(i,u)$ for the random walk $Y_N^*(t)$. At this point, it is enough to observe that $\varphi_{N,D_N}^*$ equals the r.h.s. of (3.7).

Now fix $\varepsilon > 0$. Then, due to (3.3) and (3.7),

(3.8)
$$\begin{aligned}|\Pi_N^{(1)}(t,t_w) &- \Pi_N(t,t_w)| \\ &\leq \mathbb{P}_N(x_N(t_w) \geq \varepsilon) \\ &\quad + \sum_{j:\, x_j < \varepsilon} \mathbb{P}_N(Y_N(t_w) = j) \frac{x_j}{N} \sum_{i \in D_N} \int_0^t \varphi_{N,D_N}(i, t-s)\, ds \\ &\leq \mathbb{P}_N(x_N(t_w) \geq \varepsilon) + \varepsilon \int_0^t e^{-\delta u(1 - |D_N|/N)}\, du.\end{aligned}$$

By the law of large numbers,

$$\lim_{N \uparrow \infty} \int_0^\infty e^{-\delta u(1 - |D_N|/N)} < \infty \qquad \text{a.s.}$$

The proposition now follows from the fact that $\varepsilon$ is arbitrary and from (3.5).

To deal with case (ii), one proceeds in essentially the same way, decomposing the path of the process at its returns to the point $x_i$, and summing over the number of these returns. One finds easily that the case when the process does not leave $x_i$ for the entire period $t$ dominates, leading to the assertion of the proposition. We leave the details to the reader. □

**4. The REM-like trap model on a Poisson point process.** In this section we consider a slightly different formulation of the REM like trap model that betrays more directly its connection to the REM dynamics (see [3, 4]) and that offers a somewhat more natural insight in the role of time-scales in the analysis of aging systems. Let us consider a Poisson point process, $\mathcal{P} = \sum_i \delta_{E_i}$, on $\mathbb{R}$ with intensity measure $\alpha e^{-\alpha E}\, dE$, where $0 < \alpha < 1$. Note that such processes arise naturally as the extremal process of sequences of random variables. Almost surely, the support of $\mathcal{P}$ is an infinite set of points, whose maximal element is finite. Thus, we can label the points in the support of $\mathcal{P}$ in decreasing order: $E_1 > E_2 > \cdots$. The energy landscape, $\underline{E}$, is defined as $\underline{E} = (E_1, E_2, \ldots)$. We want to define a random process on the support of this point process that jumps "uniformly" from any point to any other point in the support. To do this, we need to introduce a cut-off. Here we fix an energy threshold $E$ and set

$$N_E = \max\{i : E_i \geq E\}.$$

Note that $N_E$ is a Poisson random variable with expectation $e^{-\alpha E}$. Moreover, the probability that $N_E = 0$ can be made as small as desired when $E$ is chosen small enough, as we assume in what follows.

Let $\mathcal{G}_E = (\mathcal{S}_E, \mathcal{E}_E)$ be the graph with

$$\mathcal{S}_E := \{1, 2, \ldots, N_E\}, \qquad \mathcal{E}_E := \{\{i,j\} : i \neq j \in \mathcal{S}_E\}.$$

Since we want to investigate the effect of *time rescaling*, we introduce a time unit, $\tau_0 = e^{E_0}$. Then, the continuous-time random walk, $Y_E(t)$, is the random walk on $\mathcal{G}_E$ having uniform initial distribution and such that, after arriving at site $i \in \mathcal{E}_E$, it waits an exponential time with mean $\frac{N_E}{N_E-1} e^{E_i}/\tau_0$ and then jumps with uniform probability to a different site of $\mathcal{S}_E$. In particular, the Markov generator $\mathbb{L}_E$ for the above defined random walk is given by $\mathbb{L}_N$ in (2.3) with $N := N_E$ and $x_i := \tau_0 e^{-E_i}$ (since for $E \ll 0$, $\frac{N_E}{N_E-1} \sim 1$ when referring to waiting time, we disregard the coefficient $\frac{N_E}{N_E-1}$, as in Section 2). Note that $Y_E(t)$ depends on $\tau_0$, but we do not make this explicit in our notation. In what follows we denote by $\mathbb{P}_E$ the probability measure on the path space determined by $Y_E(\cdot)$, and by $\mathbb{E}_E$ the corresponding expectation.

Note that the physical waiting time (the absolute one) for the system at state $i$ is given by $T_i := e^{E_i}$, while in the above dynamics the waiting time is $\tau_i := T_i/\tau_0$, in agreement with the choice to consider $\tau_0$ as our new time unit. In what follows we consider, when taking the thermodynamic limit $E \downarrow -\infty$, three different kinds of time rescaling: $\tau_0$ fixed, $\tau_0 := e^E$ (i.e., $E_0 = E$) and $\tau_0 \downarrow 0$ after $E \downarrow -\infty$.

As in Section 2, we are interested in the asymptotic behavior of time–time correlation functions. In particular, let us introduce here the correlation function

$$\Pi_E(t, t_w) := \mathbb{P}_E(Y_E(s) = Y_E(t_w), \ \forall s \in [t_w, t_w + t]).$$

We will prove that, when $\tau_0$ is fixed, the system exhibits fast relaxation, thus excluding aging behavior (see Proposition 4.2). At the other extreme, the scaling $\tau_0 = e^E$ corresponds to the implicit choice made in the standard Bouchaud model considered in the previous sections. In fact, with this choice the system can be thought of as a *grand canonical* version of the original REM-like trap model and all the results of the previous sections carry over. Finally, we consider the third scaling: $\tau_0 \downarrow 0$ after $E \downarrow -\infty$. In Proposition 4.6 we show that, when performing such limits, the correlation function $\Pi_E(t, t_w)$ converges to $f(\theta)$, where $\theta = t/t_w$ and $f(\theta)$ denotes the r.h.s. of identity (2.29), that is, the limiting behavior of the correlation function $\Pi_E(t, t_w)$ is trivial. At this point, a simple consideration is fundamental. If we assume that the physical instruments in the laboratory have sensibility up to the time unit $\tau_0$, then it is natural to disregard jumps into states whose physical waiting time, $T_i = e^{E_i}$, is much smaller than $\tau_0$. Therefore, it is more appropriate to consider instead of $\Pi_E(t, t_w)$ the time–time correlation function $\Pi_E^{(1)}(t, t_w)$, defined, for $\delta > 0$ fixed, as

$$\Pi_E^{(1)}(t, t_w) = \mathbb{P}_E(x_E(u) \geq \delta \ \forall u \in (t_w, t_w + t] : x_E(u) \neq x_E(u^-)),$$

where $x_E(t) := x_k$, whenever $Y_E(t) = k$. In Section 5 we prove that $\Pi_E^{(1)}(t, t_w)$ exhibits aging behavior: $\Pi_E^{(1)}(\theta t_w, t_w)$ converges to $f(\theta)$ after taking the (ordered) limits $E \downarrow -\infty$, $\tau_0 \downarrow 0$ and $t_w \uparrow \infty$.

Finally, we discuss the asymptotic spectral behavior for the above time rescalings. We will show that aging appears whenever the limiting spectral density has a singularity of order $O(x^{\alpha-1})$ at 0.

Let us recall some properties of the Ppp $\sum_i \delta_{x_i}$ with intensity measure $\alpha \tau_0^{-\alpha} x^{\alpha-1} dx$ on $(0, \infty)$, which will be frequently used below. Given $M > 0$, the truncated Ppp, $\sum_{x_i \leq M} \delta_{x_i}$, can be realized as follows: Let $n_M$ be a Poisson variable with expectation $(\frac{M}{\tau_0})^\alpha = \int_0^M \alpha \tau_0^{-\alpha} x^{\alpha-1} dx$ and let $X_i$, $i \in \mathbb{N}$, be i.i.d. random variables on $[0, M]$ with probability distribution $p(X) dX = \alpha M^{-\alpha} X^{\alpha-1} dX$. Then

$$\text{(4.1)} \qquad \sum_{x_i \leq M} \delta_{x_i} \sim \sum_{i=1}^{n_M} \delta_{X_i},$$

in the sense that the point processes above have the same distribution. In particular, taking $M = \tau_0 e^{-E}$, we get

$$\text{(4.2)} \qquad \sum_{i \leq N_E} \delta_{x_i} \sim \sum_{i=1}^{n_E^*} \delta_{X_i},$$

where $n_E^*$ is a Poisson variable with expectation $e^{-\alpha E}$, and $X_i, i \in \mathbb{N}$, are i.i.d. random variables, independent of $n_E^*$, distributed on $[0, \tau_0 e^{-E}]$ with probability distribution $p(X) dX = e^{\alpha E} \alpha \tau_0^{-\alpha} X^{\alpha-1} dX$.

NOTATION. It is convenient to introduce the random walks $x_E(t)$, $\tau_E(t)$, defined as

$$x_E(t) := x_k, \qquad \tau_E(t) := \tau_k \qquad \text{if } Y_E(t) = k.$$

We denote by $\gamma_E$ the positive oriented loop having support

$$\text{supp}(\gamma_E) = \{x \pm i : x \in [-1, \tau_0 e^{-E} + 1]\}$$
$$\cup \{-1 + bi : |b| \leq 1\} \cup \{\tau_0 e^{-E} + 1 + bi : |b| \leq 1\}.$$

Moreover, we call $\gamma_\infty$ the infinite open path, oriented from $\infty + i$ to $\infty - i$, having support

$$\text{supp}(\gamma_\infty) = \{x \pm i : x \geq -1\} \cup \{-1 + bi : |b| \leq 1\}.$$

Finally, for given $E$, $0 = \lambda_1^{(E)} < \lambda_2^{(E)} < \cdots < \lambda_{N_E}^{(E)}$ are the $N_E$ distinct eigenvalues of the infinitesimal generator $\mathbb{L}_E$ (see Proposition 2.1).

4.1. *$\tau_0$ fixed.* Let us first observe that $\sum_{i=1}^{\infty} \tau_i < \infty$, for almost all $\underline{E}$. In fact, since the Ppp $\sum_i \delta_{\tau_i}$ has intensity measure $\alpha \tau_0^\alpha \tau^{-(1+\alpha)} d\tau$ on $(0, \infty)$,

$$\mathbb{E}(|\{i : \tau_i \geq 1\}|) = \int_1^\infty \alpha \tau_0^\alpha \tau^{-(1+\alpha)} d\tau < \infty,$$

$$\mathbb{E}\left(\sum_{i : \tau_i < 1} \tau_i\right) = \int_0^1 \alpha \tau_0^\alpha \tau^{-\alpha} d\tau < \infty.$$

Whenever $\sum_{i=1}^{\infty} \tau_i < \infty$, it is easy to derive the asymptotic spectral behavior of the system from Proposition 2.1 and to show its fast relaxation, thus implying the absence of aging:

PROPOSITION 4.1. *For almost all $\underline{E}$,*

(4.3) $$\lim_{E \downarrow -\infty} \sum_{j=1}^{N_E} \delta_{\lambda_j^{(E)}} = \sum_{j=1}^{\infty} \delta_{\lambda_j} \qquad \text{vaguely in } \mathcal{M}([0, \infty)),$$

*where $\mathcal{M}([0, \infty))$ denotes the space of locally bounded measure on $[0, \infty)$ and*

$$\{0 = \lambda_1 < \lambda_2 < \lambda_3 < \cdots\} = \left\{\lambda \in \mathbb{C} : \sum_{k=1}^\infty \frac{\lambda}{\lambda - x_k} = 0\right\}.$$

PROOF. In what follows we assume that $\sum_{i=1}^{\infty} \tau_i < \infty$, which is true a.s. Then the function $\phi_\infty(\lambda) := \sum_{k=1}^\infty \frac{\lambda}{x_k - \lambda}$ is well defined on $\mathbb{C} \setminus \{x_i : i \geq 1\}$ and has nonnegative zeros $0 = \lambda_1 < \lambda_2 < \cdots$, such that $x_{i-1} < \lambda_i < x_i$ for any $i > 1$. At this point it is enough to show that

$$\lim_{E \downarrow -\infty} \lambda_i^{(E)} = \lambda_i \qquad \forall i = 1, 2, \ldots.$$

The assertion is trivial for $i = 1$. Suppose that $N_E \geq i > 1$ and set $\psi_E(\lambda) := \sum_{k=1}^{N_E} \frac{1}{x_k - \lambda}$. Due to Proposition 2.1, $\lambda_i^{(E)}$ is the unique zero of $\psi_E(\lambda)$ in the interval $(x_{i-1}, x_i)$. In particular,

$$\psi_E(\lambda_i) = \psi_E(\lambda_i) - \psi_E(\lambda_i^{(E)}) = \int_{\lambda_i^{(E)}}^{\lambda_i} \dot{\psi}_E(\lambda) \, d\lambda.$$

Since $\dot{\psi}_E(\lambda) \geq \frac{1}{(x_i - x_{i-1})^2}$ for all $\lambda \in (x_{i-1}, x_i)$, we get

$$|\lambda_i^{(E)} - \lambda_i| \leq (x_i - x_{i-1})^2 |\psi_E(\lambda_i)|$$

and, therefore, the assertion follows by observing that the identity $\phi_\infty(\lambda_i) = 0$ implies

$$|\psi_E(\lambda_i)| \leq \sum_{k=N_E+1}^\infty \frac{1}{x_k - x_i} \downarrow 0 \qquad \text{as } E \downarrow -\infty. \qquad \square$$

PROPOSITION 4.2. *For almost all $\underline{E}$,*

$$(4.4) \quad \lim_{t\uparrow\infty}\lim_{E\downarrow-\infty}\mathbb{P}_E(x_E(t)=x_j)=\frac{\tau_j}{\sum_{k=1}^{\infty}\tau_k} \quad \forall j=1,2,\ldots,$$

*thus implying*

$$(4.5) \quad \lim_{t_w\uparrow\infty}\lim_{E\downarrow-\infty}\Pi_E(\theta t_w, t_w)=0 \quad \forall \theta>0,$$

$$(4.6) \quad \lim_{t_w\uparrow\infty}\lim_{E\downarrow-\infty}\Pi_E(t, t_w)=\frac{\sum_{i=1}^{\infty}\tau_i e^{-x_i t}}{\sum_{i=1}^{\infty}\tau_i} \quad \forall t>0.$$

PROOF. In what follows we assume that $\underline{E}$ satisfies $\sum_i \tau_i < \infty$. Setting $h(x)=\mathbb{I}_{x=x_j}$ in (2.13), we get the integral representation

$$(4.7) \quad \mathbb{P}_E(x_E(t)=x_j)=\frac{1}{2\pi i}\int_{\gamma_E}\frac{e^{-\lambda t}}{\lambda(x_j-\lambda)}\left(\sum_{k=1}^{N_E}\frac{1}{x_k-\lambda}\right)^{-1} d\lambda.$$

Applying the residue theorem [see the arguments used in order to derive (2.19)], it is easy to compute the Laplace transform, $\hat{F}_E(\omega)=\int_0^{\infty}\mathbb{P}_E(x_E(t)=x_j)e^{-\omega t}\,dt$ for $\Re(\omega)>0$:

$$(4.8) \quad \hat{F}_E(\omega)=\left(\omega(\omega+x_j)\sum_{k=1}^{N_E}\frac{1}{\omega+x_k}\right)^{-1}.$$

Using that, for $\omega=a+ib$ and $N$ is any positive integer,

$$\left|\sum_{j=1}^{N}\frac{1}{\omega+x_j}\right| \geq \begin{cases} \dfrac{1}{|\omega+x_1|}, & \text{if } a\geq 0, \\ \dfrac{b}{(a+x_1)^2+b^2}, & \text{if } a<0, \end{cases}$$

we obtain that, almost surely, there exists $c>0$, such that

$$(4.9) \quad |\hat{F}_E(\omega)|\leq c\frac{1}{|\omega|} \quad \forall E, \forall \omega\in\mathcal{A}:=\{re^{i\theta}:0<r<\infty, |\theta|\leq \tfrac{3}{4}\pi\}.$$

Let us now introduce the path, $\tilde{\gamma}$, consisting of the parabolic arcs $\{-t\pm it^2:t\geq 1\}$ and the circular arc of radius $\sqrt{2}$ around the origin connecting (in anti-clockwise way) $-1-i$ to $-1+i$. The orientation of $\tilde{\gamma}$ is such that $-1+i$ comes before $-1+i$. Then, by means of (4.9), the Laplace inversion formula and Lebesgue's dominated convergence theorem, we get

$$\lim_{E\downarrow-\infty}\mathbb{P}_E(x_E(t)=x_j)=\frac{1}{2\pi i}\int_{\tilde{\gamma}}e^{t\omega}\hat{F}(\omega)\,d\omega,$$

$$\hat{F}(\omega):=\left(\omega(\omega+x_j)\sum_{k=1}^{\infty}\frac{1}{\omega+x_k}\right)^{-1}.$$

Note that $\hat{F}(\omega)$ is the limit of $\hat{F}_E(\omega)$, as $E \downarrow \infty$; in particular, it satisfies (4.9). Moreover, $\hat{F}(\omega)$ is the Laplace transform of $\lim_{E \downarrow -\infty} \mathbb{P}_E(x_E(t) = x_j)$ and

$$\left| \omega \hat{F}(\omega) - \frac{\tau_j}{\sum_{k=1}^{\infty} \tau_k} \right| \leq c|\omega| \qquad \forall |\omega| \leq 1 : \omega \in \mathcal{A}.$$

At this point (4.4) follows from Proposition A.1. Moreover, from (4.4) and the identity

$$\Pi_E(t, t_w) = \sum_{j=1}^{N_E} \mathbb{P}_E(x_E(t_w) = x_j) e^{-(N_E - 1)/N_E x_j t},$$

one infers (4.5) and (4.6). □

4.2. $\tau_0 = e^E$. Note that, choosing $\tau_0 = e^E$, the random variables $X_1, X_2, \ldots$ introduced in (4.2) are i.i.d. with distribution given by $p(X) \, dX = \alpha X^{\alpha-1} \, dX$ on $[0, 1]$. Therefore, due to (4.2), we can think of $Y_E(t)$ as the *grand canonical* version of Bouchaud's REM-like trap model. In particular, it exhibits the same asymptotic spectral density and the same aging behavior:

PROPOSITION 4.3. *For almost all $\underline{E}$,*

$$\lim_{E \downarrow -\infty} \frac{1}{N_E} \sum_{j=1}^{N_E} \delta_{\lambda_j^{(E)}} = \alpha x^{\alpha-1} \, dx \qquad \text{weakly in } \mathcal{M}([0,1]).$$

PROOF. Approximating continuous functions on $[0, 1]$ by step functions having rational values and jumps at rational points, it is enough to prove that, given $0 \leq a < b \leq 1$,

$$\lim_{E \downarrow -\infty} \frac{1}{N_E} |\{j : 1 \leq j \leq N_E, \lambda_j^{(E)} \in [a, b]\}| = b^\alpha - a^\alpha \qquad \text{a.s.}$$

We set

$$A_E := |\{j : 1 \leq j \leq N_E, x_j \in [a, b]\}|$$
$$= |\{j : j \geq 1, e^{-E_j} \in [e^{-E}a, e^{-E}b]\}|.$$

Then, due to Proposition 2.1, we only need to prove that

$$\lim_{E \downarrow -\infty} \frac{A_E}{N_E} = b^\alpha - a^\alpha \qquad \text{a.s.}$$

To do so, observe that $N_E$ and $A_E$ are Poisson variables with parameter $e^{-\alpha E}$ and $e^{-\alpha E}(b^\alpha - a^\alpha)$. For $n \in \mathbb{N}$, we set $E(n) = -\frac{2}{\alpha} \ln n$, that is, $e^{-\alpha E(n)} = n^2$. The Chebyshev inequality and Borel–Cantelli lemma then imply that

$$\lim_{n \uparrow \infty} \frac{N_{E(n)}}{e^{-\alpha E(n)}} = 1, \qquad \lim_{n \uparrow \infty} \frac{A_{E(n)}}{e^{-\alpha E(n)}} = b^\alpha - a^\alpha \qquad \text{a.s.}$$

By monotonicity, one can extend the first limit to $\lim_{E \downarrow -\infty} \frac{N_E}{e^{-\alpha E}} = 1$, a.s. In order to extend the second limit to general $E$, we observe that, whenever $E(n+1) < E \leq E(n)$,

$$|A_E - A_{E(n)}| \leq |\{j : e^{-E_j} \in [ae^{-E(n)}, ae^{-E(n+1)}] \cup [be^{-E(n)}, be^{-E(n+1)}]\}|.$$

Since the r.h.s. is a Poisson variable with expectation of order $O(n)$, the Chebyshev inequality and the Borel–Cantelli lemma imply

$$\lim_{n \uparrow \infty} \sup_{E(n+1) < E \leq E(n)} \frac{|A_E - A_{E(n)}|}{e^{-\alpha E(n)}} = 0 \qquad \text{a.s.},$$

which allows, to prove that $\lim_{E \downarrow -\infty} \frac{A_E}{e^{-\alpha E}} = b^\alpha - a^\alpha$, a.s. □

PROPOSITION 4.4. *For almost all $\underline{E}$, and $\gamma$ a positive oriented loop around $[0, 1]$,*

$$\lim_{E \downarrow -\infty} \Pi_E(t, t_w)$$

(4.10)
$$= \frac{1}{2\pi i} \int_\gamma \frac{e^{-t_w \lambda}}{\lambda} \frac{\int_0^1 (e^{-xt}/(\lambda - x)) x^{\alpha-1} \, dx}{\int_0^1 (1/(\lambda - x)) x^{\alpha-1} \, dx} \, d\lambda \qquad \forall t, t_w.$$

*In particular, for almost all $\underline{E}$, given $\theta > 0$,*

$$(4.11) \quad \lim_{t_w \uparrow \infty} \lim_{E \downarrow -\infty} \Pi_E(\theta t_w, t_w) = \frac{\sin(\pi \alpha)}{\pi} \int_{\theta/(1+\theta)}^1 u^{-\alpha}(1-u)^{\alpha-1} \, du.$$

PROOF. Our starting point is (4.2) and the following inequality, which holds for any bounded function, $f$, with $\mathbb{E}(f(X_i)) = 0$:

$$\mathbb{P}(|\mathrm{Av}_{j=1}^k f(X_j)| \geq \delta) \leq 2 \exp\left(-\frac{k\delta^2}{4\|f\|_\infty}\right) \qquad \forall \delta > 0, k = 1, 2, \ldots.$$

In particular, conditioning on $n_E^*$ [see (4.2)], we get

$$(4.12) \qquad \mathbb{P}(|\mathrm{Av}_{j=1}^{n_E^*} f(X_j)| \geq \delta) \leq 2 \exp\{-e^{-\alpha E}(1 - e^{-\delta^2/4\|f\|_\infty})\}.$$

(4.10) can now be derived from (4.12), the Borel–Cantelli lemma, and the integral representation

$$(4.13) \qquad \Pi_E(t, t_w) = \frac{1}{2\pi i} \int_\gamma \frac{e^{-t_w \lambda}}{\lambda} \frac{\mathrm{Av}_{j=1}^{N_E} e^{-x_j t}/(x_j - \lambda)}{\mathrm{Av}_{j=1}^{N_E} 1/(x_j - \lambda)} d\lambda,$$

where $\gamma$ is a positive oriented closed path around $[0, 1]$ [see (2.12)]. Note that the r.h.s. of (4.10) corresponds to the function $\Pi(t, t_w)$ introduced in Proposition 2.5. Therefore, the assertion of Proposition 4.4 follows from Propositions 2.5 and 2.8. □

4.3. $\tau_0 \downarrow 0$ *after* $E \downarrow -\infty$. In this scaling regime, we show that the vague limit of the suitably rescaled spectral density is given by the measure $\alpha x^{\alpha-1} \, dx$ on $[0, \infty)$ and we recover the aging property of the correlation function. Moreover, due to the fact that we are effectively already at "infinite times" on the microscopic scale, we get a pure aging function even before taking $t$ and $t_w$ to infinity:

PROPOSITION 4.5. *For almost all $\underline{E}$,*

$$\lim_{\tau_0 \downarrow 0} \lim_{E \downarrow -\infty} \tau_0^\alpha \sum_{j=1}^{N_E} \delta_{\lambda_j^{(E)}} = \alpha x^{\alpha-1} \, dx \qquad \text{vaguely in } \mathcal{M}([0, \infty)).$$

PROPOSITION 4.6. *For almost all energy landscapes, $\underline{E}$, given positive $t, t_w$,*

$$(4.14) \qquad \lim_{E \downarrow -\infty} \Pi_E(t, t_w) = \frac{1}{2\pi i} \int_{\gamma_\infty} \frac{e^{-t_w \lambda}}{\lambda} \frac{\tau_0^\alpha \sum_{j=1}^\infty e^{-x_j t}/(x_j - \lambda)}{\tau_0^\alpha \sum_{j=1}^\infty 1/(x_j - \lambda)} \, d\lambda$$

*and*

$$(4.15) \quad \lim_{\tau_0 \downarrow 0} \lim_{E \downarrow -\infty} \Pi_E(t, t_w) = \frac{\sin(\pi \alpha)}{\pi} \int_{\theta/(1+\theta)}^1 u^{-\alpha}(1-u)^{\alpha-1} \, du \qquad \text{where } \theta = \frac{t}{t_w}.$$

REMARK. The integral in (4.14) exists due to Lemma 4.8.

Due to Proposition 2.1, Proposition 4.5 follows, if one is able to prove that $\tau_0^\alpha \sum_{j=1}^{N_E} \delta_{x_j}$ converges vaguely to $\alpha x^{\alpha-1} \, dx$ on $[0, \infty)$ when taking the (ordered) limits $E \downarrow -\infty$, $\tau_0 \downarrow 0$. This is the content of Lemma 4.7 below (which is analogous to Lemma 4.16 in [4]). Finally, the proof of Proposition 4.6 is based on (and given after) the technical Lemmas 4.7 and 4.8.

LEMMA 4.7. *Let $M > 0$ and let $f$ be a bounded, continuous function on $[0, M]$. Then, there exists $\delta > 0$, such that, for almost all energy landscapes, $\underline{E}$,*

$$(4.16) \qquad \left| \tau_0^\alpha \sum_{x_i \leq M} f(x_i) - \int_0^M f(x) \alpha x^{\alpha-1} \, dx \right| \leq c \tau_0^\delta \qquad \forall \tau_0 > 0,$$

*where $c > 0$ is a positive constant.*

PROOF. Let $X_1, X_2, \ldots$ and $n_M$ be as in (4.1). Due to (4.1) and since $\text{Var}(n_M) = \mathbb{E}(n_M) = (M/\tau_0)^\alpha$,

$$\mathbb{P}\left( \left| \frac{|\{j : x_j \leq M\}|}{(M/\tau_0)^\alpha} - 1 \right| \geq \varepsilon \right) \leq (\tau_0/M)^\alpha \varepsilon^{-2}.$$

In particular, given $\gamma, s > 0$ such that $2s - \gamma\alpha < -1$, using the Borel–Cantelli lemma, we obtain that, for almost all energy landscapes, there exists $c > 0$, such that
$$\left|\frac{|\{j : x_j \leq M\}|}{(M/\tau_0)^\alpha} - 1\right| \leq ck^{-s} \qquad \forall k = 1, 2, \ldots, \text{where } \tau_0 := k^{-\gamma}.$$

Due to this estimate, we get that

(4.17) $$\left|\tau_0^\alpha \sum_{x_i \leq M} f(x_i) - M^\alpha \mathrm{Av}_{x_i \leq M} f(x_i)\right| \leq ck^{-s}\|f\|_\infty$$
$$\forall k \in 1, 2, \ldots, \text{where } \tau_0 := k^{-\gamma},$$

where $\mathrm{Av}_{x_i \leq M}$ denotes the average over the set $\{x_i \leq M\}$. As done for (4.12), if $0 < \rho < 1$,

$$\mathbb{P}\left(\left|\mathrm{Av}_{x_i \leq M} f(x_i) - M^{-\alpha} \int_0^M f(x)\alpha x^{\alpha-1}\, dx\right| \geq \rho\right)$$
$$\leq 2\exp\left\{-\left(\frac{M}{\tau_0}\right)^\alpha (1 - e^{-c\rho^2})\right\}$$
$$\leq 2e^{-c'\rho^2 \tau_0^{-\alpha}}.$$

In particular, by the Borel–Cantelli lemma, for almost all $\underline{E}$,

(4.18) $$\left|\mathrm{Av}_{x_i \leq M} f(x_i) - M^{-\alpha} \int_0^M f(x)\alpha x^{\alpha-1}\, dx\right| \leq ck^{-s}$$
$$\forall k = 1, 2, \ldots, \text{where } \tau_0 := k^{-\gamma},$$

if $s$ is chosen small enough. Now (4.17) and (4.18) imply the assertion of the lemma, if $\tau_0 = k^{-\gamma}$, for some $k = 1, 2, \ldots$. The general case $\tau_0 > 0$ follows easily from the uniform continuity of $f$. □

LEMMA 4.8. *For almost all energy landscapes, $\underline{E}$, there are positive constants, $\tau_0^*, c_1, c_2$, such that the following holds: If $\tau_0 \leq \tau_0^*$, $N \geq |\{j : x_j \leq 1\}|$ and $\lambda \in \gamma_\infty$ (or $\lambda = a + ib$, with $|b| \leq 1$ and $a \geq x_j + 1$, for all $j \leq N$), then*

(4.19) $$\left|\tau_0^\alpha \sum_{j=1}^N \frac{1}{x_j - \lambda}\right| \geq c_1|\lambda|^{-2}.$$

*Moreover, if $\tau_0 \leq \tau_0^*$ and $M \geq 1$, then*

(4.20) $\tau_0^\alpha \sum_{j=1}^\infty \dfrac{1}{|x_j - \lambda|} \leq c_2|\lambda|^{\alpha-1}\ln(1+|\lambda|) \qquad$ *if $\lambda \in \gamma_\infty$,*

(4.21) $\tau_0^\alpha \sum_{x_j \geq M} \dfrac{1}{|x_j - \lambda|} \leq c_2 M^{\alpha-1}\ln M \qquad$ *if $\lambda \in \gamma_\infty$ or $\Re(\lambda) = M + 1$,*

where $\sum_{x_j \geq M}$ means $\sum_{j \geq 1 \,:\, x_j \geq M}$.

PROOF. It is convenient to introduce the nonrescaled Ppp $\sum_i \delta_{y_i}$, where $y_i := e^{-E_i}$, with intensity measure $\alpha y^{\alpha-1} \, dy$ on $(0, \infty)$. We set $x_j = \tau_0 y_j$. Moreover, we fix $\beta > 2$ and $0 < \gamma < \beta/2 - 1$ and we define

$$N_n := |\{j : n^{\beta/\alpha} \leq y_j < (n+1)^{\beta/\alpha}\}|,$$

for $n$ positive integer. Then the Borel–Cantelli lemma implies that, for almost all $\underline{E}$,

$$(4.22) \qquad \left|\frac{N_n}{(n+1)^\beta - n^\beta} - 1\right| \leq c n^{-\gamma} \qquad \forall\, n = 1, 2, \ldots.$$

Moreover, again using the Borel–Cantelli lemma and a simple argument based on monotonicity, one can prove that there exists $\delta > 0$, such that, for almost all $\underline{E}$,

$$(4.23) \qquad \left|\frac{|\{j : y_j \leq u\}|}{u^\alpha} - 1\right| \leq \kappa u^{-\delta} \qquad \forall\, u \geq 1,$$

where $\kappa$ is a positive constant. We leave this proof as an exercise.

In what follows we write $\lambda = a + ib$. Then

$$(4.24) \qquad \left|\sum_{j=1}^N \frac{1}{x_j - \lambda}\right|^2 = \left(\sum_{j=1}^N \frac{x_j - a}{(x_j - a)^2 + b^2}\right)^2 + \left(\sum_{j=1}^N \frac{b}{(x_j - a)^2 + b^2}\right)^2.$$

In order to prove (4.19), we assume (4.22) and (4.23) to be valid and we let $0 < \tau_0 \leq \tau_0^* \leq 1$, where $\tau_0^*$ is such that $1 > \kappa(\tau_0^*)^\delta$. This implies that $\{x_j : x_j \leq 1\} \neq \varnothing$. By (4.24), if $N \geq |\{j : x_j \leq 1\}|$ and $\lambda \in \gamma_\infty$,

$$\tau_0^\alpha \left|\sum_{j=1}^N \frac{1}{x_j - \lambda}\right|$$

$$\geq \begin{cases} \tau_0^\alpha \sum_{x_j \leq 1} \frac{|b|}{(x_j - a)^2 + b^2} \geq c|\lambda|^{-2} \tau_0^\alpha |\{x_j : x_j \leq 1\}|, & \text{if } |b| \geq \tfrac{1}{2}, \\ \tau_0^\alpha \sum_{x_j \leq 1} \frac{x_j + 1}{(x_j + 1)^2 + b^2} \geq c \tau_0^\alpha |\{x_j : x_j \leq 1\}|, & \text{if } |b| < \tfrac{1}{2}. \end{cases}$$

At this point, (4.19), for $\lambda \in \gamma_\infty$, follows from (4.23). The case $\lambda = a + ib$, with $|b| \leq 1$ and $a \geq x_j + 1$, for all $j \leq N$, can be treated similarly.

It is easy to derive (4.20) and (4.21) from the estimates (4.25)–(4.31) below that hold for almost all $\underline{E}$:

if $a \leq 100$, $1 \leq M$ and $\lambda = a + ib \in \gamma_\infty$, then

$$\tau_0^\alpha \sum_{j=1}^\infty \frac{1}{|x_j - \lambda|} \leq c, \tag{4.25}$$

$$\tau_0^\alpha \sum_{x_j \geq M} \frac{1}{|x_j - \lambda|} \leq cM^{\alpha-1}; \tag{4.26}$$

if $\tau_0$ is small enough and $a \geq 100$, then

$$\tau_0^\alpha \sum_{x_j \leq a/2} \frac{1}{|x_j - a|} \leq ca^{\alpha-1}, \tag{4.27}$$

$$\tau_0^\alpha \sum_{a/2 \leq x_j \leq a-1} \frac{1}{|x_j - a|} \leq ca^{\alpha-1} \ln a, \tag{4.28}$$

$$\tau_0^\alpha \sum_{a-1 \leq x_j \leq a+1} \frac{1}{|x_j - \lambda|} \leq ca^{\alpha-1} \qquad \text{if } \lambda = \alpha + ib \in \gamma_\infty, \tag{4.29}$$

$$\tau_0^\alpha \sum_{a+1 \leq x_j \leq 2a} \frac{1}{|x_j - a|} \leq ca^{\alpha-1} \ln a, \tag{4.30}$$

$$\tau_0^\alpha \sum_{x_j \geq T} \frac{1}{|x_j - a|} \leq cT^{\alpha-1} \qquad \text{if } T \geq 2a. \tag{4.31}$$

Let $\lambda \in \gamma_\infty$ with $a \leq 100$. Then, due to (4.23),

$$\tau_0^\alpha \sum_{x_j \leq 1} \frac{1}{|x_j - \lambda|} \leq c\tau_0^\alpha \left| \left\{ j : y_j \leq \frac{1}{\tau_0} \right\} \right| \leq c',$$

while, due to (4.22),

$$\tau_0^\alpha \sum_{x_j \geq 1} \frac{1}{|x_j - \lambda|} \leq c\tau_0^\alpha \sum_{x_j \geq 1} \frac{1}{x_j} \leq c'\tau_0^{\alpha-1} \sum_{n \geq \lfloor \tau_0^{-\alpha/\beta} \rfloor} n^{\beta-1-\beta/\alpha} \leq c'', \tag{4.32}$$

which proves (4.25). The proof of (4.26) follows from the same arguments. (4.27) is a simple consequence of (4.23). The l.h.s. of (4.29) can be bounded by $\tau_0^\alpha |\{j : a - 1 \leq \tau_0 y_j \leq a + 1\}|$ and (4.22) allows to conclude the proof of (4.29).

The proof of (4.28), (4.30) and (4.31) can be easily derived from the following estimate. Let $1 \leq A \leq B$ with $B \leq a - 1$ or $A \geq a + 1$, then (4.22) implies

$$\tau_0^\alpha \sum_{A \leq x_i \leq B} \frac{1}{|x_i - a|} \leq c\tau_0^\alpha \sum_{n=n_-}^{n=n_+} \frac{n^{\beta-1}}{|\tau_0 n^{\beta/\alpha} - a|} \leq c'\tau_0^\alpha \int_u^v \frac{x^{\beta-1}}{|a - \tau_0 x^{\beta/\alpha}|} dx$$

$$= c'a^{-1+\alpha} \int_{u(\tau_0/a)^{\alpha/\beta}}^{v(\tau_0/a)^{\alpha/\beta}} \frac{y^{\beta-1}}{|1-y^{\beta/\alpha}|}\, dy,$$

where $n_- = \lfloor (A/\tau_0)^{\alpha/\beta} \rfloor - 1$, $n_+ = \lfloor (B/\tau_0)^{\alpha/\beta} \rfloor + 1$, $u = n_- - 1$, $v = n_+ + 1$ (we assume $\tau_0$ small enough in order to exclude the singular point in the above intervals of sum and integral).  □

PROOF OF PROPOSITION 4.6. In order to avoid confusion, we underline here the dependence on $\tau_0 = e^{E_0}$ by writing $\Pi_{E,E_0}(t, t_w)$ instead of $\Pi_E(t, t_w)$. Our starting point is given by the integral representation (2.12):

$$(4.33) \qquad \Pi_{E,E_0}(t, t_w) = \frac{1}{2\pi i} \int_{\gamma_E} \frac{e^{-t_w \lambda}}{\lambda} \frac{\tau_0^\alpha \sum_{j=1}^{N_E} e^{-x_j t}/(x_j - \lambda)}{\tau_0^\alpha \sum_{j=1}^{N_E} 1/(x_j - \lambda)}\, d\lambda.$$

Let us choose $\underline{E}$ satisfying Lemma 4.8. Then, due to the exponential decaying factor $e^{-t_w \lambda}$ and Lemma 4.8, if $\tau_0 \leq \tau_0^*$, and $E$ is small enough, such that $\tau_0 e^{-E} \geq 1$, the integration path $\gamma_E$ in (4.33) can be replaced by $\gamma_\infty$. At this point, (4.14) follows from Lemma 4.8 and Lebesgue's dominated convergence theorem.

To prove (4.15), given a positive integer $M$, we set

$$g_{M,E_0}(t, t_w) := \frac{1}{2\pi i} \int_{\Gamma_M} \frac{e^{-t_w \lambda}}{\lambda} \frac{\tau_0^\alpha \sum_{x_j \leq M} e^{-x_j t}/(x_j - \lambda)}{\tau_0^\alpha \sum_{x_j \leq M} 1/(x_j - \lambda)}\, d\lambda,$$

where $\Gamma_M$ is the positive oriented path whose support is

$$\mathrm{supp}(\Gamma_M) = \{\lambda \in \mathbb{C} : |\lambda - x| = 1 \text{ for some } x \in [0, M]\}.$$

Then, applying Lemma 4.8, whenever $\tau_0 \leq \tau_0^*$,

$$\left| \lim_{E \downarrow -\infty} \Pi_{E,E_0}(t, t_w) - g_{M,E_0}(t, t_w) \right| \leq cM^{\alpha-1} \ln M \qquad \forall M \in \mathbb{N}_+.$$

Let us assume that $\underline{E}$ satisfies (4.16), for all $M \in \mathbb{N}_+$ and for $f(x) = \frac{e^{-xt}}{x-\lambda}$, or $f(x) = \frac{1}{x-\lambda}$, for all $\lambda$ in a countable dense set of $\Gamma_M$ and for all rational positive $t$. Then, by a chaining argument, we get

$$\lim_{E_0 \downarrow -\infty} g_{M,E_0}(t, t_w) = g_M(t, t_w) \qquad \forall t, t_w > 0,$$

where

$$(4.34) \qquad g_M(t, t_w) = \frac{1}{2\pi i} \int_{\Gamma_M} \frac{e^{-t_w \lambda}}{\lambda} \frac{\int_0^M e^{-xt}/(\lambda - x) x^{\alpha-1}\, dx}{\int_0^M 1/(\lambda - x) x^{\alpha-1}\, dx}\, d\lambda,$$

which implies

$$\limsup_{E_0 \downarrow -\infty} \left| \lim_{E \downarrow -\infty} \Pi_{E,E_0}(t, t_w) - g_M(t, t_w) \right| \leq cM^{\alpha-1} \ln M \qquad \forall M \in \mathbb{N}_+.$$

At this point, we observe (see the proof of Lemma 4.8) that there exist $c, c' > 0$ such that for all $M \in \mathbb{N}_+$:

$$\left| \int_0^M \frac{x^{\alpha-1}}{\lambda - x} dx \right| \geq c|\lambda|^{-2} \qquad \forall \lambda \in \Gamma_M \cup \gamma_\infty,$$

(4.35)
$$\int_0^M \frac{x^{\alpha-1}}{|\lambda - x|} dx \leq c' \qquad \forall \lambda \in \Gamma_M \cup \gamma_\infty,$$

$$\int_M^\infty \frac{x^{\alpha-1}}{|\lambda - x|} dx \leq c' M^{\alpha-1} \ln M \qquad \forall \lambda \in \gamma_\infty.$$

From the above estimates, we infer

(4.36) $$|g_M(t, t_w) - g(t, t_w)| \leq c M^{\alpha-1} \ln M,$$

where

(4.37) $$g(t, t_w) := \frac{1}{2\pi i} \int_{\gamma_\infty} \frac{e^{-t_w \lambda}}{\lambda} \frac{\int_0^\infty e^{-xt}/(\lambda - x) x^{\alpha-1} dx}{\int_0^\infty 1/(\lambda - x) x^{\alpha-1} dx} d\lambda.$$

Using the analytic properties of the integrand in the r.h.s. of (4.37), one can show that $g(t, t_w) = g(t/t_w, 1)$. In order to compute $g(\theta, 1)$, we observe that for a suitable positive constant $c$, $|g(\theta s, s) - g_M(\theta s, s)| \leq c M^{\alpha-1} \ln M$, for any $s \geq 1$ [in fact, the constant $c$ in (4.36) can be chosen uniformly if $t_w \geq 1$]. By the results of Section 2.1 [cf. (4.34) with $\Pi(t, t_w)$ in Proposition 2.5], we get

$$\lim_{s \uparrow \infty} g_M(\theta s, s) = \text{ r.h.s. of (4.15)},$$

thus concluding the proof. □

**5. Other correlation function when $\tau_0 \downarrow 0$ after $E \downarrow -\infty$.** As stated in Proposition 4.6, the standard time–time correlation function $\Pi_E(t, t_w)$ has trivial behavior after taking the limits $E \downarrow -\infty$, $\tau_0 \downarrow 0$. For physical reasons, it is more natural to disregard jumps into states with physical waiting time $T_i = e^{E_i}$ much smaller than $\tau_0$, since we assume that the physical instruments in the laboratory have sensibility up to the time unit $\tau_0$. Therefore, let us fix $\delta > 0$ and consider here the more natural time–time correlation function

$$\Pi_E^{(1)}(t, t_w) = \mathbb{P}_E(x_E(u) \geq \delta \ \forall u \in (t_w, t_w + t] : x_E(u) \neq x_E(u^-)).$$

The main result of this section is the following:

PROPOSITION 5.1. *For almost all $\underline{E}$,*

(5.1) $$\lim_{t_w \uparrow \infty} \sup_{t > 0} \limsup_{E_0 \downarrow -\infty} \limsup_{E \downarrow -\infty} |\Pi_E^{(1)}(t, t_w) - \Pi_E(t, t_w)| = 0.$$

*In particular, for almost all $\underline{E}$,*

$$\lim_{t_w \uparrow \infty} \lim_{E_0 \downarrow -\infty} \lim_{E \downarrow -\infty} \Pi_E^{(1)}(\theta t_w, t_w)$$

(5.2)
$$= \frac{\sin(\pi\alpha)}{\pi} \int_{\theta/(1+\theta)}^1 u^{-\alpha}(1-u)^{\alpha-1}\, du \qquad \forall\, \theta > 0.$$

Note that the correlation function $\Pi_E^{(1)}(t, t_w)$ is the analog of $\Pi_N^{(1)}(t, t_w)$ of Proposition 3.1. As for the proof of Proposition 3.1, a useful observation is that, given $\delta > 0$,

(5.3) $$\lim_{t \uparrow \infty} \lim_{\tau_0 \downarrow 0} \lim_{E \downarrow -\infty} \mathbb{P}_E(x_E(t) > \delta) = 0 \qquad \text{a.s.}$$

We can prove a stronger result concerning the phenomenon that with high probability the system visits deeper and deeper traps. In fact, note that, by (2.13),

$$\mathbb{P}_E(x_E(t) > \delta) = \frac{1}{2\pi i} \int_{\gamma_E} \frac{e^{-t\lambda}}{\lambda} \frac{\sum_{j=1}^{N_E} \mathbb{I}_{x_j \geq \delta}/(x_j - \lambda)}{\sum_{j=1}^{N_E} 1/(x_j - \lambda)}\, d\lambda.$$

Then, reasoning as in the proof of Proposition 4.6 and using the results of Section 2.2, one can easily show the analog of Proposition 2.9:

PROPOSITION 5.2. *For almost all $\underline{E}$,*

(5.4) $$\lim_{t \uparrow \infty} \lim_{\tau_0 \downarrow 0} \lim_{E \downarrow -\infty} t^{1-\alpha} \mathbb{P}_E(x_E(t) > \delta) = \frac{B(\delta)}{c(\alpha)},$$

*where*

$$B(\delta) := \frac{\int_\delta^\infty x^{\alpha-2}\, dx}{\int_0^\infty x^{\alpha-1}/(1+x)\, dx}, \qquad c(\alpha) := \int_0^\infty y^{\alpha-1} e^{-y}\, dy.$$

PROOF OF PROPOSITION 5.1. Trivially, (5.2) follows from (5.1) and Proposition 4.6. Since $E_0 \downarrow -\infty$ after $E \downarrow -\infty$, we assume that $E < E_0$ and define

$$D_{E, E_0} := \{i : E \leq E_i \leq E_0\}.$$

This set corresponds to the small traps into which we allow the particle to jump. By (5.3) and the same arguments used in the proof of Proposition 3.1 for $i = 1$ (with exclusion of the last step since here $\lim_{E \downarrow -\infty} \frac{|D_{E, E_0}|}{N_E} = 1$ a.s.), it is easy to derive the assertion of the proposition from Lemma 5.3. $\square$

LEMMA 5.3. *For almost all $\underline{E}$, there exist positive constants, $p, c$, independent of $E, E_0$, satisfying the following property. Whenever $|\{i: E_i > E_0\}| > 0$,*

$$\text{(5.5)} \qquad \limsup_{E \downarrow -\infty} \frac{1}{N_E} \sum_{i \in D_{E,E_0}} \varphi_{E,E_0}(i,t) \leq c e^{-pt} \qquad \forall t > 0,$$

*where $D_{E,E_0} := \{i: E \leq E_i \leq E_0\}$, for $E < E_0$, and the function $\varphi_{E,E_0}$ is defined as*

$$\text{(5.6)} \qquad \varphi_{E,E_0}(i,t) := \mathbb{P}_E(Y_E(u) \in D_{E,E_0} \; \forall u \in [0,t] | Y_E(0) = i).$$

PROOF. Let us assume that $E < E_0$, $|\{i: E_i > E_0\}| > 0$ and, without loss of generality, $\delta = 1$.

We fix $\ell > 0$ such that $e^{-\alpha \ell} < \frac{1}{2}$ and define

$$W_{1,E} := \{i: E \leq E_i < E + \ell\}, \qquad N_{1,E} := |W_{1,E}|,$$
$$W_{2,E} := \{i: E + \ell \leq E_i \leq E_0\}, \qquad N_{2,E} := |W_{2,E}|.$$

Note that $D_{E,E_0} = W_{1,E} \cup W_{2,E}$ and $N_E, N_{1,E}, N_{2,E}$ are Poisson variables with parameters $e^{-\alpha E}$, $e^{-\alpha E}(1 - e^{-\alpha \ell})$ and $e^{-\alpha E - \alpha \ell} - e^{-\alpha E_0}$, respectively. In particular, for almost all $\underline{E}$,

$$\text{(5.7)} \qquad \lim_{E \downarrow -\infty} \frac{N_E}{e^{-\alpha E}} = 1, \qquad p_1 := \lim_{E \downarrow -\infty} \frac{N_{1,E}}{N_E} = 1 - e^{-\alpha \ell},$$
$$p_2 := \lim_{E \downarrow -\infty} \frac{N_{2,E}}{N_E} = e^{-\alpha \ell} < \frac{1}{2}.$$

We observe that $\tilde{n} := N_E - N_{1,E} - N_{2,E}$ is a positive integer, independent of $E$ and $x_i \geq e^{E_0 - E - \ell}$, if $i \in W_{1,E}$, while $x_i \geq 1$, if $i \in W_{2,E}$. Let us introduce a new random walk, $Y_E^*(t)$, on $\mathcal{S}_E$, whose infinitesimal generator, $\mathbb{L}_E^*$, is defined as $\mathbb{L}_E$, with $x_i$ replaced by $x_i^*$ defined as

$$x_i^* = \begin{cases} x_i, & \text{if } i \notin D_{E,E_0}, \\ A_E := e^{E_0 - E - \ell}, & \text{if } i \in W_{1,E}, \\ 1, & \text{if } i \in W_{2,E}. \end{cases}$$

We denote by $\mathbb{P}_E^*$ the probability on path space associated to $Y_E^*(t)$ with uniform initial distribution and we set

$$\text{(5.8)} \qquad \varphi_{E,E_0}^*(i,t) := \mathbb{P}_E^*(Y_E^*(u) \in D_{E,E_0} \; \forall u \in [0,t] | Y_E^*(0) = i).$$

By a simple coupling argument, one gets $\varphi_{E,E_0}(i,t) \leq \varphi_{E,E_0}^*(i,t)$. In particular,

$$\text{(5.9)} \qquad \frac{1}{N_E} \sum_{i \in D_{E,E_0}} \varphi_{E,E_0}(i,t) \leq \Phi := \frac{1}{N_E} \sum_{i \in D_{E,E_0}} \varphi_{E,E_0}^*(i,t)$$

$$\forall i \in D_{E,E_0}, \; \forall t \geq 0.$$

At this point, it remains to estimate $\Phi$. In order to simplify notation, we write $D, N, N_1, N_2, A$, dropping the index $E$. Moreover, we consider the following realization of the dynamics of $Y_E^*$: after arriving at a site $i$, the walk waits an exponential time of parameter $x_i^*$ and then it jumps to a point of $\mathcal{S}_E$ with uniform probability. In particular, jumps can be *degenerate*, that is, initial and final sites can coincide.

We claim that
$$\Phi = \Phi_1 + \Phi_2 + \Phi_3,$$
where

(5.10)
$$\Phi_1 := \sum_{k_1=1}^{\infty} \sum_{k_2=0}^{\infty} \binom{k_1+k_2}{k_1} \left(\frac{N_1}{N}\right)^{k_1} \left(\frac{N_2}{N}\right)^{k_2+1} A^{k_1}$$
$$\times \int_0^t du\, e^{-Au} \frac{u^{k_1-1}}{(k_1-1)!} e^{-(t-u)} \frac{(t-u)^{k_2}}{k_2!},$$

(5.11)
$$\Phi_2 := \sum_{k_1=0}^{\infty} \sum_{k_2=1}^{\infty} \binom{k_1+k_2}{k_1} \left(\frac{N_1}{N}\right)^{k_1+1} \left(\frac{N_2}{N}\right)^{k_2}$$
$$\times A^{k_1} \int_0^t du\, e^{-Au} \frac{u^{k_1}}{k_1!} e^{-(t-u)} \frac{(t-u)^{k_2-1}}{(k_2-1)!},$$

(5.12) $\quad \Phi_3 := \dfrac{N_1}{N} \exp\left\{-At\left(1-\dfrac{N_1}{N}\right)\right\} + \dfrac{N_2}{N} \exp\left\{-t\left(1-\dfrac{N_2}{N}\right)\right\}.$

The above identities can be derived from the probabilistic interpretation of $k_1$, $k_2$ as
$$k_i = |\{\text{jumps performed before time } t \text{ having starting point in } W_i\}|$$
and from the following simple identities:
$$\mathbb{P}(T_1 + T_2 + \cdots + T_n \in [z, z+dz))$$
$$= e^{-\kappa z} \kappa^n \frac{z^{n-1}}{(n-1)!} dz,$$
$$\mathbb{P}(T_1 + T_2 + \cdots + T_n \leq z \text{ and } T_1 + T_2 + \cdots + T_n + T_{n+1} > z)$$
$$= e^{-\kappa z} \kappa^n \frac{z^n}{n!},$$
where $z \geq 0$ and $T_1, T_2, \ldots$ are independent exponential variables with parameter $\kappa$.

Finally, we only need to prove that, for suitable positive constant $c, p > 0$,

(5.13) $$\limsup_{E \downarrow -\infty} \Phi_i \leq c e^{-pt} \qquad \forall i = 1, 2, 3.$$

We give the proof in the case $i = 1$; the case $i = 2$ is completely similar, while the case $i = 3$ follows directly from (5.7).

We fix $\gamma : \alpha < \gamma < 1$, set $k_0 := A^\gamma$ and write

$$\Phi_1 = \Phi_1^{\leq k_0} + \Phi_1^{> k_0},$$

where $\Phi_1^{\leq k_0}$ is the contribution to $\Phi_1$ of addenda in the r.h.s. of (5.10) with $1 \leq k_1 \leq k_0$ and $k_2 \geq 0$.

If $k_1 > k_0$, then

$$\left(\frac{N_1}{N}\right)^{k_1} = \left(\frac{N - N_2}{N}\right)^{k_1} \left(1 - \frac{\tilde{n}}{N - N_2}\right)^{k_1} \leq \left(\frac{N - N_2}{N}\right)^{k_1} \left(1 - \frac{1}{N}\right)^{A^\gamma},$$

thus implying that

$$\Phi_1^{> k_0} \leq \Phi_1\left(\frac{N - N_2}{N}, N_2\right)\left(1 - \frac{1}{N}\right)^{A^\gamma} \leq \left(1 - \frac{1}{N}\right)^{A^\gamma} \downarrow 0 \qquad \text{as } E \downarrow -\infty,$$

where $\Phi_1(\frac{N-N_2}{N}, N_2)$ is defined as in the r.h.s. of (5.10) with $N_1$ replaced by $N - N_2$. Note that it does not exceed 1 since it corresponds to the probability of a certain event.

Let us now consider the term $\Phi_1^{\leq k_0}$. To this aim, since $\binom{k_1+k_2}{k_1} \leq 2^{k_1+k_2}$,

$$\Phi_1^{\leq k_0} \leq \sum_{k_1=1}^{k_0} \sum_{k_2=0}^{\infty} \left(\frac{2AN_1}{N}\right)^{k_1} \left(\frac{2tN_2}{N}\right)^{k_2} \frac{I(0,t)}{(k_1-1)! k_2!},$$

where

$$I(w_1, w_2) := \int_{w_1}^{w_2} e^{-(t-u) - Au} u^{k_1 - 1} \, du.$$

Fix $0 < m < 1$ with $2p_2 + 2mp_1 < 1$ (recall that $2p_2 < 1$). Then, trivially,

(5.14) $$I\left(0, \frac{mt}{A}\right) \leq \frac{1}{k_1} e^{-t(1-m/A)} \left(\frac{mt}{A}\right)^{k_1}.$$

From such a bound, one gets immediately

(5.15) $$\begin{aligned}\sum_{k_1=1}^{k_0} \sum_{k_2=0}^{\infty} &\left(\frac{2AN_1}{N}\right)^{k_1} \left(\frac{2tN_2}{N}\right)^{k_2} \frac{1}{(k_1-1)! k_2!} I\left(0, \frac{mt}{A}\right) \\ &\leq c \exp\left\{-t\left(1 - \frac{m}{A} - 2\frac{N_1}{N} m - 2\frac{N_2}{N}\right)\right\}.\end{aligned}$$

The last expression, when $E \downarrow -\infty$, converges to $c \exp(-t(1 - 2p_1 m - 2p_2))$, in agreement with (5.13).

In order to estimate the integral $I(\frac{mt}{A}, t) = e^{-t} \int_{mt/A}^{t} e^{-(A-1)u} u^{k_1-1} \, du$, we observe that

$$\int_s^w e^{-zu} u^n \, du = (-1)^n \frac{d^n}{dz^n} \frac{e^{-zs} - e^{-zw}}{z}$$

$$\leq \frac{e^{-zs}}{z} \left(s + \frac{1}{z}\right)^n + \frac{e^{-zw}}{z} \left(w + \frac{1}{z}\right)^n \qquad \forall s, w, z \geq 0,$$

thus implying the bound

$$I\left(\frac{mt}{A}, t\right) \leq c e^{-t} A^{-k_1} (mt+1)^{k_1-1} + e^{-At} (A-1)^{-1} \left(t + \frac{1}{A-1}\right)^{k_1-1}.$$

The contribution of $ce^{-t} A^{-k_1} (mt+1)^{k_1-1}$ to $\Phi_1^{>k_0}$ can be treated by means of estimates similar to the ones leading to (5.15).

In order to conclude, we only need to show that

$$\text{(5.16)} \qquad e^{-At} \sum_{k_1=1}^{A^\gamma} \sum_{k_2=1}^{\infty} \left(\frac{2AtN_1}{N}\right)^{k_1-1} \left(\frac{2tN_2}{N}\right)^{k_2} \frac{1}{(k_1-1)! k_2!} \downarrow 0$$

as $E \downarrow -\infty$.

To this aim, observe that

$$\text{r.h.s. of (5.16)} \leq e^{-At + 2tN_2/N} \sum_{k_1=0}^{A^\gamma} \left(\frac{2AtN_1}{N}\right)^{k_1}$$

$$\leq c(t) e^{-At} A^\gamma (4tA)^{A^\gamma} = c(t) \exp\{-At + \gamma \ln A + A^\gamma \ln(4tA)\}.$$

Since $0 < \gamma < 1$, we get (5.16), thus concluding the proof. □

## APPENDIX

### A.1. Laplace transform.

PROPOSITION A.1. *Let $G(t)$ be a bounded measurable function on $(0, \infty)$ and let us consider the Laplace transform*

$$\hat{G}(\omega) = \int_0^\infty G(t) e^{-t\omega} \, dt$$

*well defined if $\Re(\omega) > 0$. Let us define*

$$\text{(A.1)} \qquad \mathcal{A} := \{r e^{i\theta} : 0 < r < \infty, |\theta| \leq \tfrac{3}{4}\pi\}.$$

*Suppose that $\hat{G}$ can be analytically continued to $\mathbb{C} \setminus (-\infty, 0]$ and that there are positive constants $\gamma, \beta, \alpha, c$ and $B \in \mathbb{R}$ such that*

$$\text{(A.2)} \qquad |\hat{G}(\omega)| \leq c|\omega|^{-\gamma} \qquad \forall \omega \in \mathcal{A}, |\omega| \geq 1,$$

$$\text{(A.3)} \qquad |\omega^\beta \hat{G}(\omega) - B| \leq c|\omega|^\alpha \qquad \forall \omega \in \mathcal{A}, |\omega| \leq 1.$$

*Then,*

$$\lim_{s\uparrow\infty} s^{1-\beta}G(s) = \frac{B}{c(\beta)} \qquad \text{where } c(\beta) := \int_0^\infty y^{\beta-1}e^{-y}\,dy.$$

PROOF. If we set $H(s) := \frac{B}{c(\beta)}s^{\beta-1}$ with $s > 0$, then the Laplace transform $\hat{H}(\omega)$ is well defined for $\Re(\omega) > 0$, $\hat{H}(\omega) = B\omega^{-\beta}$ and, trivially, $\hat{H}(\omega)$ can be analytically continued to $\mathbb{C} \setminus (-\infty, 0]$.

By the inverse formula for Laplace transform (see Chapter 4, Section 4 in [15]), we have

$$(A.4) \qquad G(s) = \lim_{K\to\infty} \frac{1}{2\pi i}\int_{x-iK}^{x+iK} e^{s\omega}\hat{G}(\omega)\,d\omega \qquad \forall s > 0, x > 0,$$

where $\omega$ runs over the vertical path connecting $x - iK$ and $x + iK$. The above formula remains true if substituting $G$ with $H$. Therefore,

$$(A.5)\quad s^{1-\beta}G(s) - \frac{B}{c(\beta)} = \lim_{K\to\infty}\frac{s^{1-\beta}}{2\pi i}\int_{x-iK}^{x+iK} e^{s\omega}(\hat{G}(\omega) - \hat{H}(\omega))\,d\omega$$

$$\forall s > 0, x > 0.$$

Let $\rho := \min(\gamma, \beta)/2$. Fix a positive number $x$ and, given $K$ and $s$, define the following paths (see Figure 2).

$\gamma_K$ is the vertical path from $x - iK$ to $x + iK$. $\gamma_{1,+}$ is the segment from $-s^{-1} + s^{-1}i$ to $-1 + i$. $\gamma_{2,+}$ is an arc from $-1 + i$ to $-K^\rho + iK$ given by the parametrization $z(t) = -t + it^{1/\rho}$ with $t \in [1, K^\rho]$. $\gamma_{3,+}$ is the horizontal segment from $-K^\rho + iK$ to $x + iK$. For $i = 1, 2, 3$, we define the path $\gamma_{i,-}$ by considering the reflection of $\gamma_{i,+}$ w.r.t. the real axis and inverting the orientation. Let $\gamma_0$ be the positive-oriented circular arc of radius $s^{-1}$ from $-s^{-1} - s^{-1}i$ to $-s^{-1} + s^{-1}i$ crossing the axis of positive real numbers.

Note that the above paths depend on $s$ and/or $K$.

Because of analyticity, the integral over $\gamma_K$ of $e^{s\omega}\hat{G}(\omega)$ is equal to the sum of the integrals over $\gamma_{3,-}, \gamma_{2,-}, \gamma_{1,-}, \gamma_0, \gamma_{1,+}, \gamma_{2,+}, \gamma_{3,+}$. The same is valid with $\hat{G}$ replaced with $\hat{H}$.

By (A.2), we have that

$$(A.6) \qquad \int_{\gamma_{3,\pm}} |e^{s\omega}\hat{G}(\omega)||d\omega| \leq ce^{sx}K^{\rho-\gamma} \downarrow 0 \qquad \text{as } K \uparrow \infty$$

and, for a suitable rational function $f$,

$$(A.7)\quad \begin{aligned} & s^{1-\beta}\int_{\gamma_{2,\pm}} |e^{s\omega}\hat{G}(\omega)||d\omega| \\ & \leq s^{1-\beta}\int_1^\infty |e^{s(-t+it^{1/\rho})}\hat{G}(-t+it^{1/\rho})(-1 + \rho^{-1}t^{(1-\rho)/\rho}i)|\,dt \\ & \leq cs^{1-\beta}\int_1^\infty e^{-st}f(t)\,dt \leq c's^{1-\beta}e^{-s/2} \downarrow 0 \qquad \text{as } s \uparrow \infty. \end{aligned}$$

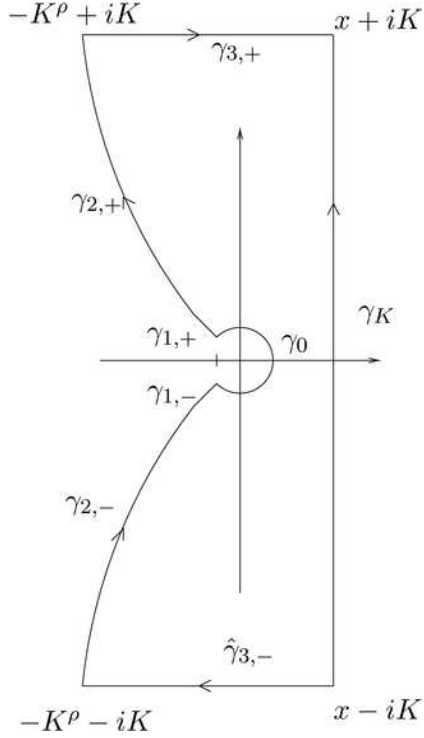

Fig. 2. *The integration path in the Laplace inversion formula.*

Similarly, it can be proved that the corresponding integrals with $\hat{G}$ substituted with $\hat{H}$ go to 0 by taking the limits $K \uparrow \infty$, $s \uparrow \infty$.

Let us now estimate $s^{1-\beta} \int_{s^{-1}}^{1} e^{-st} t^{\alpha-\beta} \, dt$ by dividing the path of integration in two paths. Choosing $0 < \delta < 1$,

$$s^{1-\beta} \int_{s^{-1}}^{s^{-1+\delta}} e^{-st} t^{\alpha-\beta} \, dt \leq c s^{1-\beta} |s^{(-1+\delta)(1+\alpha-\beta)} - s^{-(1+\alpha-\beta)}|$$

$$\leq c' s^{-\alpha} + c' s^{-\alpha+\delta(1+\alpha-\beta)},$$

$$s^{1-\beta} \int_{s^{-1+\delta}}^{1} e^{-st} t^{\alpha-\beta} \, dt \leq s^{1-\beta} e^{-s^{\delta}} g(s)$$

for a suitable rational function $g(s)$. In particular, choosing $\delta$ small enough, the above upper bounds imply

$$\lim_{s \uparrow \infty} s^{1-\beta} \int_{s^{-1}}^{1} e^{-st} t^{\alpha-\beta} \, dt = 0.$$

This result, together with (A.3), implies

$$\lim_{s \uparrow \infty} s^{1-\beta} \frac{1}{2\pi i} \int_{\gamma_{1,\pm}} e^{s\omega} (\hat{G}(\omega) - \hat{H}(\omega)) \, d\omega = 0.$$

Trivially, by (A.3),

$$s^{1-\beta}\left|\int_{\gamma_0} e^{s\omega}(\hat{G}(\omega) - \hat{H}(\omega))\, d\omega\right| \leq s^{-\alpha} \downarrow 0 \qquad \text{as } s \uparrow \infty.$$

The proposition now follows from the estimates above and (A.5). □

**A.2. Perturbation theory.** In this appendix we comment on a paper by Melin and Butaud [24] where the eigenvalues and eigenfunctions of the generator of our model were computed using perturbation theory. As pointed out earlier, these results are at variance with our exact results, and it may be worthwhile to point out the flaw in their arguments. Melin and Butaud write the generator $\mathbb{L}$ defined in (2.3) as $\mathbb{L} = T + \frac{1}{N}T^{(1)}$, where

(A.8)
$$T := \begin{pmatrix} x_1 & 0 & \cdots & 0 \\ 0 & x_2 & \cdots & 0 \\ \vdots & \vdots & \ddots & \vdots \\ 0 & 0 & \cdots & x_N \end{pmatrix},$$

$$T^{(1)} := \begin{pmatrix} -x_1 & -x_1 & \cdots & -x_1 \\ -x_2 & -x_2 & \cdots & -x_2 \\ \vdots & \vdots & \ddots & \vdots \\ -x_N & -x_N & \cdots & -x_N \end{pmatrix}.$$

The factor $1/N$ in front of the second term encourages them to consider this term as a small perturbation. Both $T$ and $T^{(1)}$ are symmetric operators on $L^2(\mu)$, where $\mu(i) := x_i^{-1}$. We denote $\langle \cdot, \cdot \rangle$ the scalar product in $L^2(\mu)$ and assume $x_1, \ldots, x_N$ to be distinct positive numbers.

Given an operator $A : L^2(\mu) \to L^2(\mu)$, we write $\|A\|$ for its operator norm. Because of symmetry, $\|T\|$ and $\|T^{(1)}\|$ are given by the maximum of $|\lambda|$, with $\lambda$ eigenvalue. Trivially, $T$ has eigenvalues $x_1, \ldots, x_N$ and $T(e_i) = x_i e_i$, where $e_1, \ldots, e_N$ is the canonical basis of $\mathbb{R}^N$, while $T^{(1)}$ has eigenvalues $0, -(x_1 + x_2 + \cdots + x_N)$.

Given $z \in \mathbb{C}$, we can define the holomorphic function $T(z) = T + zT^{(1)}$. A natural condition in order to apply perturbation theory to $T(z)$ (see [22], Chapter II) is

(A.9)
$$|z| < \frac{d}{2a_0},$$

where

$$d = \inf_{i \neq j} |x_i - x_j|, \qquad a_0 := \min_{a \in \mathbb{R}} \|T^{(1)} - a\| = \frac{x_1 + x_2 + \cdots + x_N}{2}.$$

In this case, we can conclude that $T(z) = T + zT^{(1)}$ has $N$ eigenvalues $\lambda_1(z), \ldots, \lambda_N(z)$ with $\lambda_k(z) = \sum_{n=0}^{\infty} \lambda_k^{(n)} z^n$, where, for $n \geq 2$, $|\lambda_k^{(n)}| \leq a_0^n (2/d)^{n-1}$, and

$$\lambda_k^{(1)} = \frac{\langle T^{(1)} e_k, e_k \rangle}{\langle e_k, e_k \rangle},$$

$$\lambda_k^{(2)} = \sum_{j \neq k} (x_k - x_j)^{-1} \frac{\langle T^{(1)} e_k, e_j \rangle^2}{\langle e_k, e_k \rangle \langle e_j, e_j \rangle},$$

$$\ldots.$$

Similar series exist for the perturbed eigenvectors.

However, the crucial condition (A.9) is hardly satisfied when $z = \frac{1}{N}$, since it reads

(A.10) $$\mathrm{Av}_{j=1}^{N} x_j \leq \inf_{i \neq j} |x_i - x_j|,$$

while a.s. the l.h.s. of (A.10) has nonzero limit and the r.h.s. converges to 0 at least like $1/N$.

The fact that the conditions for the application of perturbation theory are violated explains why its predictions are incorrect. This discrepancy happens not to be too obvious as far as the eigenvalues are concerned (which are caught between the diagonal elements of the generator and thus are somewhat similar to them, but the shape of the eigenfunctions is sharply different).

Namely, by Proposition 2.1, when $j \neq 1$ and $x_{j-1}, x_j$ are very near each other, the eigenvector $\psi^{(j)}$ related to the eigenvalue $\lambda_j : x_{j-1} < \lambda_j < x_j$ has two main peaks of opposite sign given by $\psi_{j-1}^{(j)}$ and $\psi_j^{(j)}$; this is very different from the predictions of [24] (see their Figure 4).

**A.3. Complex integral representation.** Let $\mathbb{L}$ be a Markov generator on the state space $\mathcal{S} := \{1, 2, \ldots, N\}$, reversible w.r.t. a positive measure $\mu$. We can think of $\mathbb{L}$ as a linear operator on $\mathbb{R}^N$, symmetric w.r.t. the scalar product $(\cdot, \cdot)_\mu$, where

$$(a, b)_\mu = \sum_{i=1}^{N} \mu(i) a_i b_i.$$

In what follows we endow $\mathbb{R}^N$ with the scalar product $(\cdot, \cdot)_\mu$ (and not with the standard Euclidean scalar product). Since $\mathbb{L}$ is symmetric, we can orthogonally decompose $\mathbb{R}^N$ as $\mathbb{R}^N = W_1 \oplus W_2 \oplus \cdots \oplus W_m$ such that $\mathbb{L} = \sum_{k=1}^{m} \lambda_k P_{W_k}$, where $P_{W_k}$ denotes the orthogonal projection on $W_k$ and $\lambda_i \neq \lambda_j$ if $i \neq j$. Given $\lambda \in \mathbb{C} \setminus \{\lambda_1, \ldots, \lambda_m\}$, we write $R(\lambda)$ for the resolvent

$$R(\lambda) := (\lambda \mathbb{I} - \mathbb{L})^{-1} = \sum_{k=1}^{m} \frac{1}{\lambda - \lambda_k} P_{W_k}.$$

Then, the residue theorem implies the integral representation

$$(A.11) \qquad e^{-t\mathbb{L}} = \frac{1}{2\pi i} \int_\gamma e^{-t\lambda} R(\lambda)\, d\lambda,$$

where $\gamma$ is a positive oriented loop containing in its interior $\lambda_1, \lambda_2, \ldots, \lambda_m$.

Given a probability measure $\nu$ on $\mathcal{S}$, we denote by $\mathbb{P}_\nu$ the probability measure on the path space associated to the continuous-time random walk $Y(t)$ on $\mathcal{S}$ with generator $\mathbb{L}$ and initial distribution $\nu$. Fix $j \in \mathcal{S}$ and let $v \in \mathbb{R}^N$ be such that $v_i = \delta_{i,j}$. We write $\frac{d\nu}{d\mu}$ for the Radon derivate, that is, $\frac{d\nu}{d\mu}(i) = \frac{\nu(i)}{\mu(i)}$. Then the symmetry of $\mathbb{L}$ w.r.t. the scalar product $(\cdot,\cdot)_\mu$ implies

$$(A.12) \qquad \begin{aligned} \mathbb{P}_\nu(Y(t)=j) &= \sum_{k=1}^N \nu(k)(e^{-t\mathbb{L}})_{k,j} \\ &= \left(\frac{d\nu}{d\mu}, e^{-t\mathbb{L}} v\right)_\mu = \left(e^{-t\mathbb{L}} \frac{d\nu}{d\mu}, v\right)_\mu \\ &= \mu(j) \sum_{k=1}^N (e^{-t\mathbb{L}})_{j,k} \frac{\nu(k)}{\mu(k)}. \end{aligned}$$

By plugging (A.11) in the r.h.s. of (A.12), we get the integral representation

$$(A.13) \qquad \mathbb{P}_\nu(Y(t)=j) = \frac{1}{2\pi i} \int_\gamma e^{-t\lambda} \left\{ \sum_{k=1}^N \mu(j) R_{jk}(\lambda) \frac{\nu(k)}{\mu(k)} \right\} d\lambda.$$

In particular, given $h$ function on $\mathcal{S}$,

$$\mathbb{E}_{\mathbb{P}_\nu}(h(Y(t))) = \sum_{j=1}^N \sum_{k=1}^N \frac{1}{2\pi i} h(j) \mu(j) \frac{\nu(k)}{\mu(k)} \int_\gamma e^{-t\lambda} R_{jk}(\lambda)\, d\lambda,$$

thus allowing to get an integral representation for $\Pi(t, t_w) := \mathbb{P}_\nu(\text{no jump in } [t_w, t_w + t])$. If we set $\mu(i) = \tau_i = x_i^{-1}$ and $\nu(i) = N^{-1}$ (uniform initial probability), then

$$(A.14) \qquad \mathbb{P}_\nu(Y(t)=j) = \frac{1}{2\pi i} \int_\gamma e^{-t\lambda} \left\{ \frac{1}{N} \sum_k \frac{x_k}{x_j} R_{jk}(\lambda) \right\} d\lambda.$$

Let us consider now the special case given by Bouchaud's REM-like trap model where $\mathbb{L} := \mathbb{L}_N$ is defined in (2.3) and $\nu$ is the uniform distribution on $\mathcal{S}$. Note that all the integral formulas obtained in Section 2 can be derived from the following one:

$$(A.15) \qquad \mathbb{P}_\nu(Y(t)=j) = \frac{1}{2\pi i} \int_\gamma e^{-t\lambda} \frac{1}{(\lambda - x_j)\phi(\lambda)}\, d\lambda,$$

where $\phi(\lambda) = \sum_{k=1}^{N} \frac{\lambda}{\lambda - x_k}$. In what follows we prove that (A.15) corresponds to (A.14).

We know already that $\det(\lambda \mathbb{I} - \mathbb{L})$ has distinct zeros given by the $N$ distinct zeros of $\phi(\lambda)$. In particular, it must be

$$(\text{A.16}) \quad \det(\lambda \mathbb{I} - \mathbb{L}) = \frac{1}{N} \phi(\lambda) \prod_j (\lambda - x_j) = \frac{1}{N} \lambda \sum_k \prod_{j:\, j \neq k} (\lambda - x_j).$$

Given a matrix $A$, we write $[A]_{i,j}$ for the determinant of the matrix obtained from $A$ by erasing the $i$th row and the $j$th column. Since

$$R_{j,k}(\lambda) = (-1)^{j+k+1} \frac{[\lambda \mathbb{I} - \mathbb{L}]_{k,j}}{\det(\lambda \mathbb{I} - \mathbb{L})}$$

and due to (A.16), in order to derive (A.15) from (A.14), we only have to show that

$$(\text{A.17}) \quad \sum_k (-1)^{j+k+1} \frac{x_k}{x_j} [\lambda \mathbb{I} - \mathbb{L}]_{k,j} = \prod_{s:\, s \neq j} (\lambda - x_s).$$

In order to prove the above identity, observe that $[\lambda \mathbb{I} - \mathbb{L}]_{k,j}$ is a polynomial of degree $N - 1$ if $k = j$, otherwise it has degree $N - 2$. The l.h.s. of (A.17) is a monic polynomial of degree $N - 1$. At this point, we only have to verify that $x_s$, $s \neq j$, are zeros of the l.h.s. of (A.17). This is trivial if one observes that the l.h.s. of (A.17) is the determinant of the matrix obtained from $\lambda \mathbb{I} - \mathbb{L}$ by replacing the $j$th column with the vector $w$ with $w_i = \frac{x_i}{x_j}$ for $i = 1, 2, \ldots, N$. It is easy to verify that, if $\lambda = x_s$ for some $s \neq j$, the $j$th row and the $s$th row in such a matrix are proportional, thus implying the thesis.

**Acknowledgments.** The authors thank J. Černý and V. Gayrard for useful discussions.

Weierstrass Institut für Angewandte
  Analysis und Stochastik
Mohrenstrasse 39
10117 Berlin
Germany
and
Mathematisches Institut
Technische Universität Berlin
Strasse des 17. Juni 136
10623 Berlin
Germany
e-mail: bovier@wias-berlin.de

Weierstrass Institut für Angewandte
  Analysis und Stochastik
Mohrenstrasse 39
10117 Berlin
Germany
e-mail: faggiona@wias-berlin.de